\documentclass[11pt]{article}
\usepackage{amsmath}
\usepackage{bbm}
\usepackage{amsfonts}
\usepackage{txfonts}
\usepackage{amssymb}
\usepackage{mathrsfs}
\usepackage{graphicx}
\usepackage{subfigure}
\usepackage{verbatim}
\usepackage{color}

\usepackage{latexsym}

\newtheorem{theorem}{\bf Theorem}[section]
\newtheorem{corollary}{\bf Corollary}[section]

\newtheorem{lemma}{\bf Lemma}[section]

\newtheorem{remark}{\bf Remark}
\newtheorem{algorithm}{Algorithm}
\numberwithin{equation}{section}
\def\udots{\mathinner{\mkern1mu\raise-1pt\vbox{\kern7pt\hbox{.}}\mkern2mu
    \raise2pt\hbox{.}\mkern2mu\raise5pt\hbox{.}\mkern1mu}}
\def\proof{{\noindent\sc Proof. \quad}}

\def\eproof{{\mbox{}\hfill\qed}\medskip}
\newcommand\qed{{\unskip\nobreak\hfil\penalty50\hskip2em\vadjust{}
\nobreak\hfil$\Box$\parfillskip=0pt\finalhyphendemerits=0\par}}
\def\proofend{\eproof}

\usepackage{algorithmic}
\usepackage[ruled]{algorithm}

\begin{document}

\title{{\bf Perturbation Analysis and Randomized Algorithms for Large-Scale Total Least Squares Problems } }
\date{}
\author{{Pengpeng Xie$^a$
\thanks{E-mail: 09110180005@fudan.edu.cn. P. Xie is supported by
the National Natural Science Foundation of China under grant 11271084.}
} \quad Yimin Wei$^{a,b}$\thanks{E-mall:
 ymwei@fudan.edu.cn. Y. Wei is supported by the National
Natural Science Foundation of China under grant 11271084.} \quad Hua
Xiang$^c$\thanks{Corresponding author (H. Xiang). E-mail:
hxiang@whu.edu.cn. H. Xiang is supported by the National Natural
Science Foundation of China under grants 10901125 and 91130022. 
Partial work was completed when he visited Fudan
University and The Chinese University of Hong Kong in 2014.
}
\\ \\
{\small $^a$ School of Mathematical Sciences, Fudan University, Shanghai, 200433, P.R. China}\\
\small{$^b$and Key Laboratory of Mathematics for Nonlinear Sciences}\\
\small{$^c$ School of Mathematics and Statistics, Wuhan University, Wuhan, 430072, P.R. China}}

\date{\today}

\maketitle

\begin{abstract}
In this paper, we present perturbation analysis and randomized
algorithms for the total least squares (\textsc{Tls}) problems. We
derive the perturbation bound and check its sharpness by numerical
experiments.
Motivated by the recently popular
probabilistic algorithms for low-rank approximations, we develop
randomized algorithms for the \textsc{Tls} and the truncated total
least squares (\textsc{Ttls}) solutions of large-scale discrete
ill-posed problems, which can greatly reduce the computational time
and still keep good accuracy.
\end{abstract}


\newpage

\section{\bf Introduction}

Given an overdetermined set of $m$ linear equations $Ax \approx b$
in $n$ unknowns $x,$ the total least squares (\textsc{Tls}) problem
can be formulated as \cite{VanHuffel}
     \begin{eqnarray}\label{eq1.1}
        \min \left\| ~[E\  f]~ \right\|_F &\textrm{subject to}& \quad b+f\in \mathscr{R}(A+E),
      \end{eqnarray}
where $\|\cdot\|_F$ denotes the Frobenius matrix norm and
$\mathscr{R} (\cdot)$ represents the range space. When the sampling or
modeling or measurement errors also affect the coefficient matrix
$A,$ the \textsc{Tls} method is more realistic, while the underlying
assumption in the least squares (\textsc{Ls}) problem is that errors
only occur in the right-hand side vector $b$.

The term ``total least squares"  was coined in \cite{Golub80}. It
has been also known as  errors-in-variables  model, orthogonal
regression, or measurement errors in the statistical literature, and
blind deconvolution in image deblurring. In the monograph \cite{VanHuffel}
the authors show the
readers how to use \textsc{Tls} for solving a variety of problems,
especially those arising in signal processing, medical imaging, and
geophysics, etc. The applications and theory associated with the
\textsc{Tls} are still being studied, for example \cite{Hansen06,
Lee}. In recent years, perturbation analysis for the \textsc{Tls}
problem has been studied extensively in the numerical linear algebra
(see e.g. \cite{Baboulin, Chang, DeMoor, Gratton2013, Hnetynkova,
Jia, Liu, Markovsky, Paige, Wei_Num, Wei, Xu, Zhou}).

It is well known that the condition number indicates the sensitivity
of the problem itself, and that an approximate bound for the forward
error can be given by the multiplication of the condition number and
the backward error. For the perturbation in the solution of the
scaled total least squares problem, Zhou et al. \cite{Zhou}
presented a first order estimation.   But as pointed out by the
authors, it is not easy to compute since the condition number
formula is a Kronecker product-based one. Baboulin and Gratton
\cite{Baboulin} derived a computable expression for the condition
number. At almost the same time, Li and Jia \cite{Jia} made a first
order perturbation analysis. Recently, Jia and Li \cite{Jia12}
proposed a formula which only used the singular values and the right
singular vectors of $[A,~b]$, and  presented the lower and upper
bounds for the condition number.  In this paper, we will present a
relative perturbation bound without considering the condition number
only.
We first give a perturbation bound in this paper.
And its significant improvements will be demonstrated by numerical examples.
We also show that these three condition numbers in \cite{Baboulin,Jia,Zhou} mentioned above are mathematically equivalent. 

For the numerical solution of the \textsc{Tls} problem,  a simple
and elegant solver based on  the \textsc{Svd} of the augmented
matrix $[A, ~b]$ can be used. When $A$ is large, a complete
\textsc{Svd} will be very costly. One improvement is to compute a
partial \textsc{Svd} based on Lanczos bi-diagonalization
\cite{Fierro}. But the partial \textsc{Svd} is still prohibitive for
large-scale sparse or structured matrices, since the initial
reduction of $[A, ~b]$ to bi-diagonal form will destroy the sparsity
or structure of the matrix. For the \textsc{Tls} problem with very
ill-conditioned coefficient matrix whose singular values decay
gradually, the task is even more challenging. Without
regularization, the ordinary least squares or total least squares
solvers yield physically meaningless solutions. For such discrete
ill-posed problems, there already exist several regularization
strategies of the \textsc{Tls} solution. For example, the solution
can be stabilized by truncating small singular values of $[A, ~b]$
via an iterative algorithm based on Lanczos bi-diagonalization
\cite{Fierro}.
Tikhonov regularization strategy is used in \cite{BeckBenTal,
GolubHansenOLeary, Lampe, Lee, LuPT}, where a Cholesky decomposition
is computed in each step in \cite{BeckBenTal}, and the linear
systems are projected onto Krylov subpace  of much smaller
$\mbox{dimensions}$ to reduce the problem size in \cite{Lampe}.
Regularization by an additional quadratic constraint is another
choice \cite{BeckBenTeboulle,LampeVoss08,Renaut,SimaHuffelGolub}, which
is the regularized \textsc{Tls} based on quadratic eigenvalue
problems (\textsc{Qep}): adding a quadratic constraint to the
\textsc{Tls}, and then iteratively solving the \textsc{Qep}.
For the large-scale discrete ill-conditioned problem, a complete
\textsc{Svd} is prohibitive, and the choice of regularization
parameter is also time consuming. The classical \textsc{Svd} of a
matrix can be well approximated by the randomized \textsc{Svd}
\cite{Halko}, and the regularization parameter can also be located
by randomized algorithms \cite{Xiang}. Such randomized algorithms
can greatly reduce the computational time, and still keep good
accuracy with very high probability. Motivated by these randomized
matrix algorithms, we present randomized algorithms for the
solution of  total least squares (\textsc{Tls}) problems
including the well-conditioned cases and
 the ill-conditioned cases.
For the practical cases where the numerical rank is not known,
randomized algorithms can usually be implemented in an adaptive approach
with the sample number increasing until the desired tolerance is satisfied.
Here the tolerance parameter is adopted to describe how well the basis matrix
generated by the randomization captures the action of the target matrix.
Based on this, we develop the randomized algorithm for the fixed precision cases.

Throughout this paper, $\mathbb{R}^{m\times n}$ denotes the set of
$m \times n$ matrices with real entries and $I_n$ stands for the
identity matrix with order $n.$ As usual, ${\bf 0}$ denotes the zero matrix
with the corresponding size easily known from the context.
For a matrix
$A\in \mathbb{R}^{m\times n},$ $A^\mathrm T$ is the transpose of
$A$; $\|A\|_2$, $\|A\|_F$ and $\|A\|_\infty$ denote the spectral norm, the
Frobenius norm and the infinity norm of $A$, respectively.
And $A^{\dagger}$ represents the Moore-Penrose inverse of $A$ \cite{Golub} and $\lambda_{\max} \left( A \right)$
denotes the largest eigenvalue of $A$.
For any matrix $A = \left[ a_1, a_2, \ldots, a_n \right]=\left( a_{ij} \right)\in \mathbb{R}^{m\times n}$
and $B = \left( b_{ij} \right) \in \mathbb{R}^{p \times q},$ the Kronecker product
 $A \otimes B$ is defined as $A \otimes B=\left( a_{ij} B \right) \in \mathbb{R}^{mp \times nq}.$
We define $\mathrm{vec}(A) = \left[ a_1^\mathrm T, a_2^\mathrm T, \ldots, a_n^\mathrm T \right]^\mathrm T \in \mathbb{R}^{mn}.$
For a vector $a,$ $\textrm{diag} (a)$ is a diagonal matrix whose diagonals are given as components of $a.$
The remaining sections of this paper are organized as follows. Section 2 introduces some basic results. In section 3, we present our
main perturbation results and show the mathematical equivalence of three kinds of condition numbers. We turn to the randomized algorithms
in section 4 and the detailed error analysis for Algorithm \textsc{Rttls} is given. The numerical results are
performed in section 5 and section 6 concludes this paper.

\section{Preliminaries}
Let $A \in \mathbb{R}^{m \times n}$ and $b \in \mathbb{R}^m$ with $m \geq n.$
Let $[A, ~b]$ and $A$  have singular value decompositions, respectively
 $$U^\mathrm T [A, ~b]V = \textrm{diag} \left( \sigma_1, \sigma_2, \ldots,\sigma_t \right)=\Sigma,$$
$$\widetilde{U}^\mathrm T A\widetilde{V}=\textrm{diag} \left( \widetilde{\sigma}_1,\widetilde{\sigma}_2,\ldots,\widetilde{\sigma}_n \right),$$
where $t = \min \{m, n+1\}$ and for the case $m > n,$ orthonormal matrices $U$ and $V,$ diagonal matrix $\Sigma$ are partitioned as follows:
       \begin{eqnarray*}
          U=\left[ U_1,u_{n+1} \right]_{m\times (n+1)}, &V=\left[\begin{array}{cc}
                          V_{11}& v_{12}\\
                          v_{21} & v_{22}
                          \end{array}\right]_{(n+1)\times(n+1)}, &\Sigma=\left[\begin{array}{cc}
                          \Sigma_1& {\bf 0}\\
                          {\bf 0}& \sigma_{n+1}
                          \end{array}\right].
      \end{eqnarray*}
For the usual well-conditioned cases in this paper, we assume the genericity condition:
\begin{equation}\label{eqn:GenericityCond}
\widetilde{\sigma}_n>\sigma_{n+1},
\end{equation}
 to ensure the
existence and uniqueness of the \textsc{Tls} solution $x$  (see
\cite{Golub80}). The singular value $\sigma_{n+1}$ can be treated
as 0 for the case $m = n$ since $\sigma_{n+1}$ does not exist at all.
From best rank-1 approximation \cite{Golub} of
matrix $[A, ~b]$, we know that
      \begin{eqnarray*}
          [E, ~f]&=&-U\left[\begin{array}{cc}
                          {\bf 0}& {\bf 0}\\
                          {\bf 0}& \sigma_{n+1}
                          \end{array}\right]V^\mathrm T\\
               &=&-\sigma_{n+1} u_{n+1} \left[v_{12}^\mathrm T, \quad v_{22}\right]\\
               &=&-\sigma_{n+1} u_{n+1} v_{n+1}^\mathrm T,
      \end{eqnarray*}
where $v_{n+1} = \left[v_{12}^\mathrm T , \quad
v_{22}\right]^\mathrm T$. Therefore, $x=-v_{12}/v_{22}.$
It follows from \cite[Theorem 2.7]{VanHuffel} that the solution $x$ can
also be expressed as a function of $[A, ~b]$, i.e.,
     \begin{eqnarray} \label{solution}
                       x & = &\left(A^\mathrm T A - \sigma_{n+1}^2 I\right)^{-1} A^\mathrm T b,
      \end{eqnarray}
and it holds that
     \begin{eqnarray} \label{x}
                       \left[x^\mathrm T, ~ -1 \right]^\mathrm T & = & -\frac{1}{v_{22}} v_{n+1}.
      \end{eqnarray}
 In the case $m = n$ under the genericity condition, the original system is just a nonsingular one
 and the \textsc{Tls} solution is equal to the least squares solution.
 Now consider the \textsc{Tls} problem (\ref{eq1.1}) and assume
 $\widetilde{\sigma}_q > \sigma_{q+1} = \cdots =\sigma_{n+1}$ with $q \leq n.$
 Let the above \textsc{Svd} still hold but partition $V$ differently as follows:
\begin{eqnarray} \label{Vpartion}
V & = &\begin{array}{@{}r@{}c@{}c@{}c@{}l@{}}
    \left[\begin{array}{c} \\ \\  \end{array}\right.
                    & \begin{array}{c} V_{11} \\ v_{21}  \end{array}
                    & \begin{array}{c}   \\    \end{array}
                    & \begin{array}{c} V_{12} \\ v_{22}  \end{array}
                          &    \left]\begin{array} {c} n \\1 \end{array}\right. \\
    & q & &n\!+\!1\!-\!q
\end{array}.
\end{eqnarray}
The condition $\widetilde{\sigma}_q > \sigma_{q+1} = \cdots =\sigma_{n+1}$ is equivalent to
that $\sigma_q > \sigma_{q+1} = \cdots =\sigma_{n+1}$ and $v_{22}$ is of full row rank,
 i.e., $v_{22}$ is not a zero vector.
According to \cite[Theorem 3.10]{VanHuffel}, the minimum norm \textsc{Tls} solution $\bar{x}$ is given by
\begin{eqnarray*}
\bar{x} &=& -V_{12}v_{22}^\dag = \left(V_{11}^\mathrm T\right)^\dag v_{21}^\mathrm T.
\end{eqnarray*}
This is called
the truncated total least squares (\textsc{Ttls}). The case $m < n+1$ requires that $\sigma_{n+1} = 0$ and
hence $\|[E, ~f]\|_F = 0$ \cite{VanHuffel}. The idea of \textsc{Ttls} is to treat the small singular values of the
augmented matrix $[A, ~b]$ as zeros and convert a numerically rank-deficient problem to an exactly rank-deficient one.
For the discrete ill-posed problems where the singular values of the coefficient matrices
decay gradually, \textsc{Ttls} can be applied, where the parameter $q$ then plays the role of the regularization parameter.
In practical applications, the smallest singular
values of $[A, ~b]$ rarely coincide \cite{Wei}. But if one considers the \textsc{Tls}
problem as an approximation to the corresponding unobservable exact relation
$A_0 x = b_0,$ then $\mbox{rank}\left([A_0, ~b_0]\right) = \mbox{rank}\left(A_0\right) = q \leq n.$
So $\sigma_{q+1}, \ldots, \sigma_{n+1}$ are just the perturbations
of zero. In this case it is realistic to define an error bound $\epsilon$ such
that all singular values $\sigma_i,$ satisfying $\left| \sigma_i - \sigma_{n+1} \right| < \epsilon,$
are considered to coincide with $\sigma_{n+1}.$ Therefore, we can use the formula
$\bar{x} = -
V_{12} v_{22}^\dag.$

\section{Perturbation results}
First, we give a lemma which will be very useful in our analysis.
        \begin{lemma}\label{lem1} Consider the total least squares problem (\ref{eq1.1}) and assume that the genericity condition \eqref{eqn:GenericityCond} holds. If $~[A, ~b]$ is perturbed to $[A + \delta A, ~b+ \delta b]$, then we have
     \begin{eqnarray*}
                      \sigma_{n+1} u_{n+1}^\mathrm T [\delta A,~\delta b] v_{n+1} & = & \frac{r^\mathrm T [\delta b - (\delta A) x]}{1+x^\mathrm T x},
      \end{eqnarray*}
where $r=b-Ax.$
        \end{lemma}
\proof
From (\ref{x}) and the singular value
decomposition of $[A,~b]$, we know that
$$
r = b-Ax = - [A,~b] \left[\begin{array}{c}
                          x\\
                          -1
                          \end{array}\right]
  = \frac{1}{v_{22}} [A,~b] v_{n+1}
  = \frac{1}{v_{22}} \sigma_{n+1} u_{n+1}.
$$
Therefore we have
     \begin{eqnarray*}
                     \frac{r^\mathrm T [\delta b - (\delta A) x]}{1+x^\mathrm T x}  & = & \frac{\sigma_{n+1}}{v_{22}}  \frac{u_{n+1}^\mathrm T [\delta b - (\delta A) x]}{1+x^\mathrm T x}\\
                     & = & - \sigma_{n+1} v_{22} u_{n+1}^\mathrm T [\delta A, ~ \delta b] \left[\begin{array}{c}
                          x\\
                          -1
                          \end{array}\right] \\
                     & = & \sigma_{n+1} u_{n+1}^\mathrm T [\delta A,~\delta b] v_{n+1},
      \end{eqnarray*}
where we use $v_{22}^2=\frac{1}{1+x^\mathrm T x}$, which is a direct result of (\ref{x}).
\proofend

The following lemma \cite{Stewart87} is also needed for deriving our
perturbation result.
        \begin{lemma}\label{lem2} Let
        $\sigma_{\min}$ be the smallest nonzero and simple singular value of a matrix $X$ with $u_{\min}$ and $v_{\min}$ being its corresponding left and right singular vectors, respectively. If $\left\| \delta X \right\|_F$ is sufficiently small, then the smallest nonzero singular value $\widehat{\sigma}_{\min}$ of the perturbed matrix $\widehat{X}=X+\delta X$ is simple
        and
     \begin{eqnarray*}
                      \widehat{\sigma}_{\min} & = & \sigma_{\min} + u_{\min}^\mathrm T (\delta X) v_{\min} + \mathcal{O}\left( \left\| \delta X \right\|_F^2 \right).
      \end{eqnarray*}

        \end{lemma}

In the following, we present our perturbation bound under the
genericity condition \eqref{eqn:GenericityCond}.
        \begin{theorem}\label{maintheorem} Consider the total least squares problem (\ref{eq1.1}) and assume that the genericity condition \eqref{eqn:GenericityCond} holds.
      If $\left\| [\delta A,~ \delta b] \right\|_F$ is sufficiently small, then we have that
     \begin{eqnarray}\label{err1}
                       \frac{\|\delta x\|_2}{\|x\|_2} & \lesssim  & \left( \frac{\|b\|_2}{\|x\|_2} \left \|\left(A^\mathrm T A - \sigma_{n+1}^2 I\right)^{-1} A^\mathrm T \right\|_2  \right) \frac{\|\delta b\|_2}{\|b\|_2} \nonumber \\
                       & + &  \left[ \frac{\|A\|_2 \|r\|_2}{\|x\|_2}  \left\| \left(A^\mathrm T A - \sigma_{n+1}^2 I\right)^{-1} \right\|_2 + \|A\|_2 \left\| \left(A^\mathrm T A - \sigma_{n+1}^2 I\right)^{-1} A^\mathrm T \right\|_2 \right] \frac{\|\delta A\|_2}{\|A\|_2},
     \end{eqnarray}
where $r=b-Ax.$
    \end{theorem}
\proof
From $\left(A^\mathrm T A - \sigma_{n+1}^2 I\right) x= A^\mathrm T
b,$ perturbing $[A,~b]$ yields
\begin{eqnarray}\label{x_pert}
\left[(A+\delta A)^\mathrm T (A+\delta A)- \widehat{\sigma}_{n+1}^2
I_n\right] (x+\delta x)= (A+\delta A)^\mathrm T (b+\delta b),
\end{eqnarray}
where $\widehat{\sigma}_{n+1}$ is the smallest singular value of
$[A+\delta A,~ b+\delta b].$ Subtracting  two equations
\eqref{solution} and \eqref{x_pert}, we have
\begin{eqnarray}\label{AAdeltax}
\left(A^\mathrm T A- \sigma_{n+1}^2 I_n\right) \delta x  & = &
(\delta A)^\mathrm T r + A^\mathrm T \left[\delta b-(\delta A) x
\right] + \left( \widehat{\sigma}_{n+1}^2 - \sigma_{n+1}^2
\right)(x+\delta x) \nonumber \\
& + & {\mathcal O}\left( \left\| [\delta A, ~ \delta b] \right\|_F^2 \right).
\end{eqnarray}
Furthermore, combining Lemma \ref{lem1} and Lemma \ref{lem2}, we have
     \begin{eqnarray} \label{sigma}
                      \widehat{\sigma}_{n+1}^2 - \sigma_{n+1}^2 & = & \left( \widehat{\sigma}_{n+1}-\sigma_{n+1}\right) \left( \widehat{\sigma}_{n+1}+\sigma_{n+1} \right) \nonumber \\
                       & = & \left\{ u_{n+1}^\mathrm T [\delta A,~ \delta b] v_{n+1} + {\mathcal O}\left( \left\|[\delta A, ~ \delta b]\right\|_F^2 \right)\right\} \left\{ 2 \sigma_{n+1} + u_{n+1}^\mathrm T [\delta A,~ \delta b] v_{n+1} + {\mathcal O}\left( \left\|[\delta A, ~ \delta b]\right\|_F^2 \right) \right\} \nonumber \\
                       & = & 2 \frac{r^\mathrm T [\delta b - (\delta A) x]}{1+x^\mathrm T x} + {\mathcal O}\left( \left\|[\delta A, ~ \delta b]\right\|_F^2 \right),
      \end{eqnarray}
where for the last approximation we use Lemma \ref{lem1}.
From (\ref{AAdeltax}) and (\ref{sigma}),
ignoring higher order terms, we know that
\begin{eqnarray*}
  \delta x
 & \approx & \left( A^\mathrm T A-\sigma_{n+1}^2 I_n \right)^{-1} \left\{ \left[ A^\mathrm T \left( \delta b - (\delta A)x \right) + (\delta A)^\mathrm T r \right] + 2 \dfrac {r^\mathrm T [\delta b-(\delta A)x]}{1+x^\mathrm T x}(x+\delta x) \right\} \nonumber\\
  & \approx & \left( A^\mathrm T A-\sigma_{n+1}^2 I_n \right)^{-1}\left[A^\mathrm T(\delta b)-A^\mathrm T(\delta A)x+(\delta A)^\mathrm T r\right] + 2 \left( A^\mathrm T A-\sigma_{n+1}^2 I_n \right)^{-1} \dfrac {r^\mathrm T \left[ \delta b-(\delta A)x \right]}{1+x^\mathrm T x}x \nonumber \\
  & = & \left( A^\mathrm T A-\sigma_{n+1}^2 I_n \right)^{-1}\left[A^\mathrm T+2\dfrac {xr^\mathrm T}{1+x^\mathrm T x}\right]\delta b
  +\left(A^\mathrm T A-\sigma_{n+1}^2 I_n\right)^{-1}\left[(\delta A)^\mathrm T r-\left (A^\mathrm T+2\dfrac {xr^\mathrm T}{1+x^\mathrm T x}\right)(\delta A)x\right].
\end{eqnarray*}
Hence,
\begin{eqnarray*}
\| \delta x \|_2
 & \lesssim & \left\| \left( A^\mathrm T A-\sigma_{n+1}^2 I_n \right)^{-1} \left[A^\mathrm T + 2\dfrac {x r^\mathrm T} {1 + x^\mathrm T x}\right] \right\|_2 \|\delta b\|_2 \\
 & + & \left\| \left( A^\mathrm T A-\sigma_{n+1}^2 I_n \right)^{-1}  \right\|_2 \|\delta A \|_2 \|r\|_2 + \left\| \left( A^\mathrm T A-\sigma_{n+1}^2 I_n \right)^{-1} \left[A^\mathrm T + 2\dfrac {x r^\mathrm T} {1 + x^\mathrm T x}\right] \right\|_2 \|\delta A\|_2 \|x\|_2.
\end{eqnarray*}
Since $r=b-Ax 
=\dfrac{\sigma_{n+1}}{v_{22}}u_{n+1}$, using the \textsc{Matlab} notation we have
\[\dfrac {xr^\mathrm T}{1+x^\mathrm T x}=-v_{n+1}(1:n) \sigma_{n+1} u_{n+1}^\mathrm T = - [V_{11},~ v_{12}]\left[\begin{array}{cc}
                          {\bf 0}& {\bf 0}\\
                          {\bf 0}& \sigma_{n+1}
                          \end{array}\right] U^\mathrm T.\]
Moreover, from the \textsc{Svd} of $[A,~ b]$, it follows that
$$A = \left [ A, ~b \right ] \left[\begin{array}{c}
                          I_n\\
                          {\bf 0}_{1\times n}
                          \end{array}\right]
     = U \Sigma \left[\begin{array}{cc}
                          V_{11}& v_{12}\\
                          v_{21} & v_{22}
                          \end{array}\right]^\mathrm T \left[\begin{array}{c}
                          I_n\\
                          {\bf 0}_{1\times n}
                          \end{array}\right] = U \Sigma [V_{11},~ v_{12}]^\mathrm T,   \qquad
                          \Sigma=\left[\begin{array}{cc}
                          \Sigma_1& {\bf 0}\\
                          {\bf 0}& \sigma_{n+1}
                          \end{array}\right].$$
Therefore, we obtain that
\begin{equation}\label{Ar}
\begin{split}
A^\mathrm T+2\dfrac {xr^\mathrm T}{1+x^\mathrm T x}
&=\left[ V_{11},~ v_{12} \right] \left[\begin{array}{cc}
                          \Sigma_1& {\bf 0}\\
                          {\bf 0}& -\sigma_{n+1}
 \end{array}\right]
U^\mathrm T .
\end{split}
\end{equation}
Furthermore,
\[\left[A^\mathrm T+2\dfrac {xr^\mathrm T}{1+x^\mathrm T x}\right]\cdot\left[A^\mathrm T+2\dfrac {xr^\mathrm T}{1+x^\mathrm T x}\right]^\mathrm T
=[V_{11},~ v_{12}] \Sigma^2 [V_{11},~ v_{12}]^\mathrm T =A^\mathrm T
A,\] and
\begin{eqnarray}\label{AArx}
\left\| \left( A^\mathrm T A-\sigma_{n+1}^2 I_n \right)^{-1} \left[ A^\mathrm T+2\dfrac {xr^\mathrm T}{1+x^\mathrm T x} \right] \right\|_2 = \left\| \left( A^\mathrm T A-\sigma_{n+1}^2 I_n \right)^{-1}A^\mathrm T \right\|_2.
 \end{eqnarray}
Finally, we have
     \begin{eqnarray*}
                       \frac{\|\delta x\|_2}{\|x\|_2} & \lesssim  & \left( \frac{\|b\|_2}{\|x\|_2} \left \|\left(A^\mathrm T A - \sigma_{n+1}^2 I\right)^{-1} A^\mathrm T \right\|_2  \right) \frac{\|\delta b\|_2}{\|b\|_2} \nonumber \\
                       & + &  \left[ \frac{\|A\|_2 \|r\|_2}{\|x\|_2}  \left\| \left(A^\mathrm T A - \sigma_{n+1}^2 I\right)^{-1} \right\|_2 + \left\| \left(A^\mathrm T A - \sigma_{n+1}^2 I\right)^{-1} A^\mathrm T \right\|_2 \|A\|_2 \right] \frac{\|\delta A\|_2}{\|A\|_2}.
     \end{eqnarray*}
\proofend

The succinct perturbation bound above is based on the formula \eqref{AArx}, which is derived
by using \eqref{Ar} and the fact that $\left\|K\right\|_2^2 = \lambda_{\max} \left( K K^\mathrm T \right)$
for any real matrix $K.$
In fact, we can give another bound of the perturbation system, and express it as the following corollary.
\begin{corollary}
Under the same conditions assumed in Theorem \ref{maintheorem}, we have
\begin{equation}\label{Corollary}
\begin{split}
\dfrac{\|\delta x\|_2}{\|x\|_2}
&\lesssim \left [\left\|\left(A^\mathrm T A-\sigma_{n+1}^2 I_n \right)^{-1}A^\mathrm T \right\|_2  \dfrac{\sqrt{ 1+\|x\|_2^2 }}{\|x\|_2} + \left\|\left(A^\mathrm T A-\sigma_{n+1}^2 I_n \right)^{-1} \right\|_2 \dfrac{\|r\|_2}{\|x\|_2} \right] \left\| \left[\delta A, ~ \delta b\right] \right\|_2. \\
\end{split}
\end{equation}
\end{corollary}
\proof
From the proof of Theorem \ref{maintheorem}, we know
\begin{equation*}\label{dx_t}
\begin{split}
\delta x  & = \left( A^\mathrm T A-\sigma_{n+1}^2 I_n \right)^{-1} \left( A^\mathrm T+2\dfrac {xr^\mathrm T}{1+x^\mathrm T x} \right)\delta b\\
  &\quad+\left(A^\mathrm T A-\sigma_{n+1}^2 I_n\right)^{-1} \left[ (\delta A)^\mathrm T r-\left( A^\mathrm T+2\dfrac {xr^\mathrm T}{1+x^\mathrm T x}\right) (\delta A)x \right] + {\mathcal O} \left( \left\|[\delta A, ~ \delta b]\right\|_F^2 \right),
\end{split}
\end{equation*}
which can be rewritten as
     \begin{eqnarray*}
                       \delta x & =  & \left( A^\mathrm T A-\sigma_{n+1}^2 I_n \right)^{-1} \left[ A^\mathrm T+2\dfrac {xr^\mathrm T}{1+x^\mathrm T x} \right] [\delta A, ~ \delta b] \left[\begin{array}{c}
                          -x\\
                          1
                          \end{array}\right]  \\
                       & +  & \left [ \left( A^\mathrm T A-\sigma_{n+1}^2 I_n \right)^{-1}, ~ {\bf 0}_{n\times 1}\right] [\delta A, ~ \delta b]^\mathrm T r  + {\mathcal O}\left( \left\|[\delta A, ~ \delta b]\right\|_F^2 \right).
     \end{eqnarray*}
%
Taking 2-norm on both sides, considering the property \eqref{AArx} and omitting the higher order terms, we simply get the bound for the relative error
\begin{equation*}
\begin{split}
\dfrac{\|\delta x\|_2}{\|x\|_2}
&\lesssim \left [\left\|\left(A^\mathrm T A-\sigma_{n+1}^2 I_n \right)^{-1}A^\mathrm T \right\|_2  \dfrac{\sqrt{ 1+\|x\|_2^2 }}{\|x\|_2} + \left\|\left(A^\mathrm T A-\sigma_{n+1}^2 I_n \right)^{-1} \right\|_2 \dfrac{\|r\|_2}{\|x\|_2} \right] \left\| \left[\delta A, ~ \delta b\right] \right\|_2.\\
\end{split}
\end{equation*}
\proofend

\begin{remark}
In Theorem \ref{maintheorem}, if $m = n,$ the relative error estimate \eqref{err1} can be simplified as
$$\frac{\|\delta x\|_2}{\|x\|_2}  \leq  \frac{\left\| A^{-1} \right\|_2 \|b\|_2}{\|x\|_2} \frac{\|\delta b\|_2}{\|b\|_2} + \|A\|_2 \left(\left\| A^{-1} \right\|_2 + \left\|A^{-1}\right\|_2^2 \frac{\|r\|_2}{\|x\|_2} \right) \frac{\|\delta A\|_2}{\|A\|_2},$$
which is one specific case in the estimate of the least squares solution \cite{Malyshev}.
From the proof of Theorem \ref{maintheorem}, we know that
\begin{eqnarray*}
\left( A^\mathrm T A-\sigma_{n+1}^2 I_n \right)^{-1} (\delta A)^\mathrm T r  & = & \left( A^\mathrm T A-\sigma_{n+1}^2 I_n \right)^{-1} A^\mathrm T  \left(A^\mathrm T\right)^\dag (\delta A)^\mathrm T r \\
& = & \left( A^\mathrm T A-\sigma_{n+1}^2 I_n \right)^{-1} A^\mathrm T \left[r^\mathrm T \otimes \left(A^\dag \right)^\mathrm T \right] \mathrm{vec}\left( \left( \delta A \right)^\mathrm T \right)
\end{eqnarray*}
and therefore the term
$\left\| \left( A^\mathrm T A-\sigma_{n+1}^2 I_n \right)^{-1} (\delta A)^\mathrm T r  \right\|_2 \leq \left\| \left( A^\mathrm T A-\sigma_{n+1}^2 I_n \right)^{-1} A^\mathrm T\right\|_2 \|r\|_2 \left\| A^\dag \right\|_2 \left\| \delta A \right\|_F.$ So we can get another bound
     \begin{eqnarray*}
                       \|\delta x\|_2 & \lesssim  & \left\| \left( A^\mathrm T A-\sigma_{n+1}^2 I_n \right)^{-1} A^\mathrm T\right\|_2 \left[ \|\delta b\|_2 + \left\| A^\dag \right\|_2  \|r\|_2 \|\delta A\|_F + \|x\|_2 \|\delta A\|_2 \right].
     \end{eqnarray*}
This bound is  succinct but it is   bigger than the bound in
\eqref{err1}.
It is easy to check
that $$\left \|\left(A^\mathrm T A - \sigma_{n+1}^2 I\right)^{-1}
A^\mathrm T \right\|_2 =
\frac{\widetilde{\sigma}_n}{\widetilde{\sigma}_n^2 -
\sigma_{n+1}^2}.$$
We notice that the term $\left( A^\mathrm T A - \sigma_{n+1}^2 I \right)^{-1} A^\mathrm T$ also appears
in the derivation of the ``effective condition number" of the total least
squares problem.  The effective condition number is defined as
\cite{Li07,Li10comp,Li10nlaa}
 $$\textrm{Cond\_eff} = \frac{\|b\|_2}{\sigma_r \|x\|_2}=  \frac{\left\| A^{\dagger} \right\|_2\|b\|_2}{\|x\|_2}$$ for the linear system $Ax=b$ with $\sigma_r$ being the smallest positive singular value of $A.$ In some cases, the effective condition number is much smaller than the traditional one.
\end{remark}

Denote
\begin{eqnarray*}
          M&=&\left[K\otimes b^\mathrm T - x^\mathrm T \otimes \left( KA^\mathrm T \right)-K\otimes (Ax)^\mathrm T,  \quad KA^\mathrm T \right],\\
          N&=& 2\sigma_{n+1} y \left(v_{n+1}^\mathrm T \otimes u_{n+1}^\mathrm T\right)
\end{eqnarray*}
with $K=\left(A^\mathrm T A-\sigma_{n+1}^2 I \right)^{-1}$ and
$y=Kx.$ Omitting the  complicated higher order term $R(\delta A,
~\delta b)$ in \cite[Eqn.(3.5)]{Zhou}, then the upper bound derived in
\cite{Zhou} becomes
$$\frac{\|M+N\|_2 \|[A, ~
b]\|_F}{\|x\|_2} \frac{\|[\delta A,~\delta b]\|_F}{\|[A, ~b]\|_F},$$
and $K_\textsc{zlwq} =\frac{\|M+N\|_2 \|[A, ~
b]\|_F}{\|x\|_2}$ can be defined as the condition number.

Denote $\widetilde{D}=\textrm{diag}
\left(\left(\widetilde{\sigma}_1^2-\sigma_{n+1}^2\right)^{-1},
\ldots,
\left(\widetilde{\sigma}_n^2-\sigma_{n+1}^2\right)^{-1}\right)$ and $
D = \textrm{diag} \left( \sqrt{\sigma_1^2 + \sigma_{n+1}^2}, \ldots,
\sqrt{\sigma_n^2 + \sigma_{n+1}^2}\right).$
The upper bound obtained from \cite{Baboulin} is expressed by $$\sqrt{ 1+\|x\|_2^2 } \left\| \widetilde{D}~[\widetilde{V}^\mathrm T, \quad {\bf 0}_{n \times 1}] ~ V ~ [D, \quad {\bf 0}_{n\times 1}]^\mathrm T \right\|_2 \frac{\|[A, ~b]\|_F}{\|x\|_2} \frac{\|[\delta A,~\delta b]\|_F}{\|[A, ~b]\|_F},$$
and they define
$$K_\textsc{bg}=\sqrt{ 1+\|x\|_2^2 } \left\| \widetilde{D}~\left[ \widetilde{V}^\mathrm T, \quad {\bf 0}_{n \times 1} \right] ~ V ~ \left[ D, \quad {\bf 0}_{n\times 1} \right]^\mathrm T \right\|_2 \frac{\|[A, ~b]\|_F}{\|x\|_2}$$ as the relative condition number.

Later, Li and Jia \cite{Jia} established the following bound for the relative perturbation  $$K_\textsc{lj} \frac{\|[\delta A, ~\delta b]\|_F}{\|[A, ~b]\|_F},$$
where $$K_\textsc{lj} = \frac{ \left \| \left( A^\mathrm T A - \sigma_{n+1}^2 I_n \right)^{-1} \left(2A^\mathrm T \frac{r}{\|r\|_2} \frac{r^\mathrm T}{\|r\|_2} G(x) - A^\mathrm T G(x) + \left [I_n\otimes r^\mathrm T , ~~ {\bf 0}_{n\times m}\right ]\right) \right\|_2 \|[A, ~
b]\|_F}{\|x\|_2}$$ is the condition number with $G(x)=\left[x^\mathrm T , ~ -1\right] \otimes I_m.$

We need to point out that, to derive the expressions for $K_\textsc{zlwq},~ K_\textsc{bg}$ and $K_\textsc{lj},$ the higher order terms have been omitted in \cite{Zhou, Baboulin, Jia}. And it is reasonable to compare our bound given in Theorem \ref{maintheorem} with the above three bounds. The numerical results will be given later.

\begin{remark}
Note that  $K_\textsc{bg}$ has another closed formula
\cite{Baboulin}
$$K_\textsc{bg} = \sqrt{ 1+\|x\|_2^2 }  \left \|\left(A^\mathrm T A - \sigma_{n+1}^2 I \right)^{-1} \left[A^\mathrm T A + \sigma_{n+1}^2 \left(I_n - \frac{2xx^\mathrm T}{1 + \|x\|_2^2}\right) \right]  \left(A^\mathrm T A - \sigma_{n+1}^2 I \right)^{-1}\right\|_2^{1/2} \frac{\|[A, ~b]\|_F}{\|x\|_2}.$$
Since $A^\mathrm T A + \sigma_{n+1}^2 \left(I_n - \frac{2xx^\mathrm
T}{1 + \|x\|_2^2}\right) =  A^\mathrm T A - \sigma_{n+1}^2 I_n + 2
\sigma_{n+1}^2 \left(I_n - \frac{xx^\mathrm T}{1 +
\|x\|_2^2}\right)$ is symmetric positive definite, we can define $L
L^\mathrm T$ as its Cholesky factorization. Then we have
$$K_\textsc{bg} = \sqrt{ 1+\|x\|_2^2 } \left \|\left(A^\mathrm T A - \sigma_{n+1}^2 I \right)^{-1} L \right\|_2 \frac{\|[A, ~b]\|_F}{\|x\|_2},$$
which is another expression of $K_\textsc{lj}$ \cite{Jia12}.

Moreover, using Lemma \ref{lem1} and the proof in \cite[Lemma 3.2]{Zhou}, we can get the following equation
      \begin{eqnarray*}
          M + N & = &\left[ -x^\mathrm T \otimes D_{\sigma_{n+1}^2} + \left(A^\mathrm T A - \sigma_{n+1}^2 I \right)^{-1} \otimes r^\mathrm T, \quad
          D_{\sigma_{n+1}^2} \right],
      \end{eqnarray*}
where $D_{\sigma_{n+1}^2}=\left(A^\mathrm T A - \sigma_{n+1}^2 I \right)^{-1} \left(A^\mathrm T + 2 \frac{x r^\mathrm T }{1+x^\mathrm T x} \right).$
Denote $P \in \mathbb{R}^{mn \times mn}$ the permutation matrix that represents the matrix transpose by
$\mathrm{vec}(B^\mathrm T) = P \mathrm{vec} (B).$ Note that $K_\textsc{bg}$ can also be expressed by \cite{Baboulin}
$$K_\textsc{bg} = \frac{\left\| \mathcal {M}_{g^\prime} \right\|_2 \|[A,~b]\|_F}{\|x\|_2}$$ with
$$\mathcal {M}_{g^\prime} = \left[ -x^\mathrm T \otimes D_{\sigma_{n+1}^2} + \left( r^\mathrm T \otimes \left(A^\mathrm T A - \sigma_{n+1}^2 I \right)^{-1} \right) P, \quad
          D_{\sigma_{n+1}^2} \right].$$
We can easily check that $M+N={\mathcal {M}}_{g^\prime},$ which means
that $K_\textsc{bg} = K_\textsc{zlwq}.$

Therefore, we see that the condition numbers derived respectively in \cite{Baboulin, Jia, Zhou} are mathematically equivalent.  But as pointed by the authors themselves, the normwise condition number proposed in \cite{Zhou} is not easy to compute.
\end{remark}

\section{Randomized algorithms for \textsc{Tls} problems}


The randomized algorithms have been receiving increasingly more attention in numerical linear algebra,
and they open the possibility of dealing with truly
massive data sets, and have become more and more popular in the
matrix approximation in the last decade \cite{Halko}. Numerical
experiments and detailed error analysis show that these random
sampling techniques can be quite effective and more efficient than the classical competitors in many
aspects. Avron et al. in \cite{Avron}  derived a randomized
least-squares solver which outperforms \textsc{Lapack} by large
factors for dense highly overdetermined systems. Recently, Xiang and
Zou \cite{Xiang} used the randomized strategy for solving
large-scale discrete inverse problems. In this section, we first propose
two algorithms for the cases where the numerical rank is known
using the similar randomized strategies. One is the
randomized algorithm for total least squares (\textsc{Rtls} for
short), and the other is the randomized algorithm for truncated
total least squares (\textsc{Rttls} for short). For the circumstances in which the
target rank is not known, we further develop the adaptive randomized $\mbox{algorithms}$ under
the fixed precision (\textsc{Arttls} for short). These randomized
algorithms can greatly reduce the computational time, and still yield
good approximate solutions.

\subsection{Randomized algorithm \textsc{Rtls} for well-conditioned cases}

\begin{algorithm} (\emph{\textsc{Algorithm Rtls: Randomized algorithm for Tls}})
\begin{itemize}
\item[1.] Generate an $(n+1) \times l$ Gaussian random matrix $\Omega
$.
\item[2.] Solve $\left( C^{\mathrm T} C \right) X = \Omega$, where $C = [ A, ~b] \in \mathbb{R}^{m \times
(n+1)}$.
\item[3.] Compute the $(n+1) \times l$ orthonormal matrix $Q$ via QR factorization $X = Q R $.
\item[4.] Solve $\left( C^{\mathrm T} C \right) Y = Q$.
\item[5.] Form the $l \times l$ matrix $Z = Q^{ \mathrm T} Y = Q^{ \mathrm T} \left(C^\mathrm T C\right)^{-1} Q$.
\item[6.] Compute the \textsc{Svd} of the smaller symmetric matrix, $Z =  W \Sigma W^{\mathrm T}$, where $W$ is orthogonal.
\item[7.] Form the $(n+1) \times l$ matrix $V = Q W$, and define
$v=V(:,1)$.
\item[8.] Form the solution $x_\textsc{rtls} = -v(1:n) /
v(n+1)$.
\end{itemize}
\end{algorithm}

For the total least squares problem, it is very crucial to obtain
the right singular vector $v_{n+1}$ associated with the smallest
singular values of $[A, ~b]$. Then the total least squares solution
can be expressed by \eqref{x}. For the expression \eqref{x}, we know
that the key point is to find the singular vector associated with
the smallest singular value. But the randomized \textsc{Svd}
\cite{Halko} usually approximates well the largest singular values
and the corresponding singular vectors. Suppose $C = [A, ~b]$ has
the \textsc{Svd} $C=U \Sigma V^\mathrm T,$ where $\Sigma =
\textrm{diag}\left( \sigma_1, \sigma_2, \ldots,\sigma_{n+1}
\right),$ $U$ and $V$ are orthogonal matrices. If $\sigma_{n+1} =
0,$ then $b$ is in the range of $A$ and the \textsc{Tls} solution is
equal to the least squares solution. We do not consider this trivial
case here. Then $C^\mathrm T C = V\Sigma^\mathrm T \Sigma V^\mathrm
T,$ and $\left(C^\mathrm T C\right)^{-1} = V \textrm{diag} \left(
\sigma_{n+1}^{-2}, \ldots,\sigma_1^{-2} \right) V^\mathrm T.$ Hence
we can see that $\sigma_{n+1}^{-2}$ becomes the largest diagonal
element, and we can apply the randomized algorithm to approximate
this value and achieve its corresponding singular vector $v_{n+1}.$
The essential step of this traditional algorithm is the \textsc{Svd}
of $[A, ~b]$. But when the size of $A$ is large, \textsc{Svd} can be
very costly, or even prohibitive. How to reduce the computational
cost and still ensure the accuracy of the approximate solution is
our main concern. Our new randomized algorithm \textsc{Rtls} is
presented in Algorithm \textsc{Rtls}.

Note that $l$ is a pre-specified parameter. In \cite{Halko} the
index $l$ is usually selected in the form $l=k+p$, where $p$ is an
oversampling parameter, and $k$ corresponds to the rank $k$
specified in advance for the best rank-$k$ approximation of $A$. To
understand Algorithm \textsc{Rtls} more, we make some remarks about
each step of the algorithm.  In Step 2 we
obtain $X = \left( C^{\mathrm T}C \right)^{-1} \Omega$ to extract the column information, 
which is further represented by an orthogonal matrix $Q$ in Step 3.
The linear system involving $C^{\mathrm T}C$ in Step 2 and 4 can be
solved by direct methods such as Cholesky factorization or Krylov
subspace iterative methods. When the problem is not too
ill-conditioned, this coefficient matrix is symmetric positive
definite, and can be solved quite efficiently. After Step 5 the
problem is reduced to a smaller symmetric semi-positive definite
matrix $Z=Q^{\mathrm T} \left( C^{\mathrm T}C \right)^{-1} Q$, and
\textsc{Svd} is applied to this small matrix in Step 6. This leads
to an \textsc{Svd} approximation, $\left( C^{\mathrm T}C
\right)^{-1} \approx V \Sigma V^{\mathrm T}$, where $V=QW$ and $W
\Sigma W^{\mathrm T} = Z$. We then use this approximate \textsc{Svd}
to seek the approximate total least squares solution
$x_\textsc{rtls}$ in Step 8.

\subsection{Randomized algorithm \textsc{Rttls} for ill-conditioned cases}

Algorithm \textsc{Rtls} works well for the well-conditioned cases.
For the total least squares problem with very ill-conditioned
coefficient matrices, the condition number of $C^\mathrm T C$ can be
very large since the condition number $\mathrm{Cond} \left( C^\mathrm T C \right) = \mathrm{Cond} \left( C \right)^2$.
We need to use regularization techniques to avoid noise
contaminations and obtain a meaningful approximate solution.
Fierro et al. in \cite{Fierro} focused on the truncated \textsc{Tls}
for solving discrete ill-posed problems, where the singular values
of the coefficient matrix decay gradually. The technique of
truncated \textsc{Tls} is similar in spirit to truncated
\textsc{Svd} (\textsc{Tsvd}), where the small singular values of
$[A, ~b]$ are treated as zeros, and the problem is reduced to an
exactly rank-deficient one \cite{Fierro}. Recently, the sensitivity
analysis and conditioning has been given in \cite{Gratton2013} and some
applications of the truncated \textsc{Tls} are reported \cite{Gratton2013}:
System identification, linear system theory, image reconstruction, speech and audio
processing, modal and spectral analysis, chemometrics, computer vision, machine learning,
computer algebra, and astronomy.
 The traditional truncated total least squares solution is given
 by the following Algorithm \textsc{Ttls} \cite[Section 3.6.1]{VanHuffel}.

\begin{algorithm} (\emph{\textsc{Algorithm Ttls:  Classical truncated Tls}})
\begin{itemize}
\item[1.] Compute the \textsc{Svd}: $[A, ~b] = U \Sigma V^\mathrm T =
\sum\limits_{i=1}^{n+1} \sigma_i u_i  v_i^\mathrm T,$ where $A \in \mathbb{R} ^ {m \times n}.$
\item[2.] Partition the matrix, $V = \begin{bmatrix} V_{11} & V_{12}
\\ v_{21} & v_{22} \end{bmatrix}$, where $V_{12} \in \mathbb{R}^{n \times
(n+1-k)},~ v_{21} \in \mathbb{R}^{1 \times k}$, and $v_{22} \in \mathbb{R}^{1 \times
(n+1-k)}$.
\item[3.] Form the minimum-norm \textsc{Tls} solution: $x_\textsc{ttls} = -
V_{12} v_{22}^\dag$.
\end{itemize}
\end{algorithm}

In Algorithm \textsc{Ttls}, the truncation parameter $k$ is
user-specified or determined adaptively \cite{Fierro}. It is chosen
such that the first $k$ large singular values dominate and
$\left\| v_{22} \right\|_2 \neq 0$. Here the Moore-Penrose inverse $v_{22}^\dag =
v_{22}^ \mathrm T ||v_{22}||_2^{-2}$.

When the discrete ill-posed problems is of medium size, we can
compute the complete \textsc{Svd} of $[A, ~b]$ directly like Step 1
in Algorithm \textsc{Ttls}. When the size of $A$ is large, the
\textsc{Svd} in Step 1 is very costly since the \textsc{Svd} needs
about $6mn^2 + 20n^3$ flops \cite{Golub}. This flaw leads us to
improve the efficiency by computing the \textsc{Svd} of $[A, ~b]$ in
Step 1 ``partially.'' The corresponding algorithm is named ``partial
total least squares (\textsc{Ptls})'' in \cite{VanHuffel}. The only
difference between \textsc{Ttls} and \textsc{Ptls} lies in the Step
1: one uses the classical complete \textsc{Svd}, while the other one
applies the partial \textsc{Svd}. The authors in \cite{VanHuffel}
report that \textsc{Ptls} is two times faster than \textsc{Ttls}
while the same accuracy can be maintained. Moreover, the relative
efficiency of partial \textsc{Svd} increases when the dimension of
the desired singular subspace is relatively smaller to the dimension
$n.$
 For large-scale  discrete
ill-posed problems, Lanczos bi-diagonalization in \cite{Fierro}  is
used to achieve a good approximation to the singular triplets
associated with several largest singular values. This approach will
lose the sparsity or structure of the coefficient matrix in the
first step of bi-diagonal reduction. What's more, Lanczos procedure
needs to access the coefficient matrix many times and use the
\textsc{Blas}-2 operations, i.e., the matrix-vector multiplications.
Here we propose an alternative technique based on randomized
strategies, that is, a randomized version of truncated total least
squares (\textsc{Rttls}). This is a new randomized algorithm, most
flops spent on the matrix-matrix multiplications, which are the
so-called nice \textsc{Blas}-3 operations, and the algorithm can be
realized by accessing the original large-scale matrix $A$ only once.

\begin{algorithm} (\emph{\textsc{Algorithm Rttls: Randomized algorithm for truncated Tls}})
\begin{itemize}
\item[1.]   Generate an $(n+1) \times l$ Gaussian random matrix $\Omega $.
\item[2.]   Form the $m \times l$ matrix $Y = C \Omega $, where $C=[A, ~b]$.
\item[3.]   Apply QR decomposition to $Y$, i.e.,  $Y = Q R $, where  $Q \in \mathbb{R}^{m \times l}$. 
\item[4.]   Form the $l \times (n+1)$ matrix $Z$ such that $Z = Q^\mathrm T C$.
\item[5.]   Apply \textsc{Svd} to the smaller matrix $Z$, i.e., $Z =  W \Sigma V^\mathrm T$, where $V \in \mathbb{R}^{(n+1)\times l}$.
\item[6.]   Let $V_{11}=V\text{(1:n, 1:$k$)}, v_{21}=V\text{(end, 1:$k$)}$, and form the solution
$x_\textsc{rttls} =  \left( V_{11}^\mathrm T \right)^\dag
v_{21}^\mathrm T $.
\end{itemize}
\end{algorithm}

Usually the randomized algorithm cannot approximate the small
singular values very well, hence we do not prefer to use the
expression $x_\textsc{ttls} = - V_{12} v_{22}^\dag$ directly. Since
in Algorithm \textsc{Rttls} we can obtain a good approximation of
the right singular vectors associated with largest singular values,
we use $x_\textsc{rttls} =  \left( V_{11}^\mathrm T \right)^\dag
v_{21}^\mathrm T $ in Step 6.
In Algorithm \textsc{Rttls} the parameter $l$ stands for the number
of sampling, and the number $k$ is the parameter for truncating ($k
\leq l$). A larger $l$ will improve the reliability of the algorithm
\cite{Halko}, but also increase the computational complexity. In
practice, we choose $l \ll n$, and make a balance between the
reliability and the computational complexity. The truncation
parameter $k$ can be user-specified or determined by some
regularization technique if no a priori estimate. Here we use
randomized regularization techniques in \cite{Xiang} to obtain an
estimation for this
parameter. 
We first perform randomized algorithms to obtain an approximate
\textsc{Svd} of $A$, then a \textsc{Gcv} function based on this
approximation is used to determine the truncation parameter $k$ for
the \textsc{Tsvd} solution of $A x \approx b$. This procedure can be
performed very fast  \cite{Xiang}.
This parameter cannot be the optimal for the total least squares
based on the \textsc{Svd} of the augmented matrix $[A, ~b]$, but
should be a reasonable estimate for the truncation parameter in
\textsc{Ttls}. Other rules such as the L-curve, quasi-optimality,
and discrepancy principle can be also used for regularization
parameter choice.
Our randomized \textsc{Ttls} is constructed in the spirit of the
truncated \textsc{Svd} (\textsc{Tsvd}). Tikhonov regularization for
\textsc{Ttls} \cite{BeckBenTal, GolubHansenOLeary, LampeVoss08, Lee,
LuPT} can be also combined with randomized algorithms, together with
the existing rules for regularization parameter choice, such as
L-curves, \textsc{Gcv}, quasi-optimality, and discrepancy principle,
etc. The detailed discussion about some important issues, such as
the regularization parameter choice, the scaling of $A$ and $b$
\cite[Section 3.6.2]{VanHuffel}, is beyond the scope of this paper.

\begin{table}[!htbp]
\begin{center}

\begin{tabular}{clllll}\hline
 Step & & \textsc{Rtls } & \textsc{Rttls} & \textsc{Ttls}  & \textsc{Ptls }\\
\hline
 1 & & ${\mathcal O}\left( nl \right)$    & ${\mathcal O}(nl)$           &   $6mn^2 + 20n^3$                & ${\mathcal O}(mnl)$            \\
 2 & & $2mn^2 + \frac{2}{3}n^3$           & $2 mnl$                      &   -                                  &   -         \\
 3 & & $4 n l^2 - \frac{4}{3} l^3$        & $4 m l^2 - \frac{4}{3} l^3$  &  ${\mathcal O}\left( n^2-nk \right)$  &  ${\mathcal O}\left( n^2 - nk \right)$   \\
 4 & & $\frac{2}{3}n^3$                   & $2 mnl$                      &   -                                   &  -          \\
 5 & & $2 n l^2$                          & $6 n l^2 + 20 l^3 $          &   -                                   & -            \\
 6 & & $ 26 l^3 $                         & ${\mathcal O}\left( n k^2 \right)$  &   -                            &   -        \\
 7 & & $2 n l^2$                          & -                            &   -                                  &  -          \\
 8 & & ${\mathcal O}(n)$                  & -                            &   -                                  &  -          \\ \hline
\end{tabular}
\caption{\label{tab:complexity_rttls}Computational complexity. }
\end{center}
\end{table}

\subsection{Adaptive randomized algorithms for truncated \textsc{Tls}}
The randomized algorithms discussed above are used to solve the fixed-rank problems.
In practical applications, the target rank is rarely known in advance.
We do not need to determine it accurately.
 So the adaptive approach \cite{Halko} is usually implemented to increase the number of samples until the
error $\left\| C - QQ^\mathrm T C \right\|_2$ satisfies the desired tolerance. The tolerance parameter
is a standard for measuring whether the basis matrix $Q$ captures the action of the target matrix $C.$
The theoretical basis behind
this scheme is that we can estimate the exact error $\left\| C - QQ^\mathrm T C \right\|_2$ by computing
$\left\| \left( I - QQ^\mathrm T \right)C \omega \right\|_2$ with $\omega$ being a standard Gaussian vector.
Draw $r$ standard Gaussian vectors, then
\begin{eqnarray*}
    \left\| \left( I - QQ^\mathrm T \right)C \right\|_2 & \leq &  10 \sqrt{\frac{2}{\pi}}  \max\limits_{1 \leq i \leq r}  \left\| \left( I - QQ^\mathrm T \right)C \omega_i \right\|_2
\end{eqnarray*}
holds except with probability $1-10^{-r}$ where $r$ is an integer that balances computational cost and reliability \cite{Halko}.

Given the augmented matrix
$C = [A, ~b],$ a tolerance $\epsilon,$ and integer $r,$ the formal schemes for computing an orthonormal basis
in Step 3 of Algorithm \textsc{Rttls} and therefore finding the truncated solution are described in Algorithm
\textsc{Arttls}. It follows that $\left\| C - QQ^\mathrm T C \right\|_2 \leq \epsilon$ holds with probability
at least $1 - \min\{m,n+1\}10^{-r}.$ We need to stress that the reorthogonalization is implemented in Step 6 and Step
7 to overcome the numerical instability that the column vectors of $Q$ become small as increasing the basis.
The CPU time requirements of Algorithms \textsc{Arttls} and \textsc{Rttls}
are essentially identical \cite{Halko}.

\begin{algorithm} (\emph{\textsc{Algorithm Arttls: Adaptive randomized algorithm for truncated Tls}})
\begin{itemize}
\item[1.] Generate standard Gaussian random vectors $\omega_1, \ldots, \omega_r$ of length $n.$
\item[2.] For $i = 1, \ldots, r,$ compute $y_i = C \omega_i.$
\item[3.] Set $j = 0$  and $Q^{(0)}= [~],$ i.e., the $m \times 0$ empty matrix.
\item[4.] while $\max \left\{\left\|y_{j+1}\right\|_2, \ldots, \left\|y_{j+r}\right\|_2 \right\} \geq \epsilon / \left( 10\sqrt{2/\pi} \right),$
\item[5.] \quad $j = j+1.$
\item[6.] \quad Overwrite $y_j$ by $\left[ I - Q^{(j-1)} \left(Q^{(j-1)}\right)^\mathrm T \right] y_j.$
\item[7.] \quad $q_j = y_j / \left\|y_j\right\|_2.$
\item[8.] \quad $Q^{(j)} = \left[ Q^{(j-1)},~~ q_j\right].$
\item[9.] \quad Draw a standard Gaussian random vector $\omega_{j+r}$ of length $n.$
\item[10.] \quad $y_{j+r} = \left[ I - Q^{(j)} \left(Q^{(j)}\right)^\mathrm T \right] C \omega_{j+r}.$
\item[11.]   \quad $\left[ y_{j+1},\ldots,y_{j+r-1} \right] = \left[ y_{j+1},\ldots,y_{j+r-1} \right] - q_j q_j^\mathrm T \left[ y_{j+1},\ldots,y_{j+r-1} \right]. $
\item[12.]   end while
\item[13.]   $Q = Q^{(j)}.$
\item[14.]   Form the $j \times (n+1)$ matrix $Z = Q^\mathrm T C$.
\item[15.]   Apply \textsc{Svd} to the smaller matrix $Z$, i.e., $Z =  W \Sigma V^\mathrm T$, where $V \in \mathbb{R}^{(n+1)\times j}$.
\item[16.]   Let $V_{11}=V\text{(1:$n$,1:$j$)}, v_{21}=V\text{(end, 1:$j$)}$, and form the solution
$x_\textsc{arttls} =  \left( V_{11}^\mathrm T \right)^\dag
v_{21}^\mathrm T $.
\end{itemize}
\end{algorithm}

\subsection{Computational complexity}

We shall say a few words about the computational complexity of the
randomized algorithms. The cost of each step of algorithms is listed
in Table \ref{tab:complexity_rttls}.  In Table
\ref{tab:complexity_rttls}, we denote the \textsc{Ptls} 
 the corresponding
 Algorithm \textsc{Ttls} using the Lanczos bi-diagonalization to fulfill
  the partial \textsc{Svd} in Step 1. As discussed in the above subsection, we
  have to form the solution by exploiting the singular vectors corresponding
  to the largest singular values, i.e.,
  $x_\textsc{ttls} = \left(V_{11}^\mathrm T\right)^\dag v_{21}^\mathrm T$.
For the matrix $C=[A, ~b] \in
\mathbb{R}^{m \times (n+1)}$, the flops count of the classical
\textsc{Svd} based on R-bidiagonalization is about $6mn^2 + 20n^3$
\cite{Golub}, while the cost of Algorithm \textsc{Rtls} is  about
$2mn^2 + \frac{4}{3}n^3 + 8n l^2 + {\mathcal O}(l^3)$; the cost of
Algorithm \textsc{Rttls} is about $4mnl + (4m+6n) l^2 + {\mathcal
O}(l^3)$. The cost of Algorithm \textsc{Rtls} is much
cheaper than the classical one. Note that the most flops are
performed in Step 2 of \textsc{Rttls} by very efficient \textsc{Blas}-3 operations,
and that fast Krylov subspace iterative solvers can be used in Step
2 and 4 of \textsc{Rtls} instead. We can see that the computational cost of
\textsc{Ptls } is of the same magnitude as \textsc{Rttls}.
But \textsc{Ptls } just carries out \textsc{Blas}-2 operations and it will be
not as efficient as it looks in practical computations.
 The advantage of our \textsc{Rtls} can be more
obvious than just what the flops account tells.

For the cases where singular values decay rapidly, we can choose a
small parameter $l$. For most cases, $m \gtrsim n \gg l$. According
to the flops, the ratio of the cost for Algorithm \textsc{Rttls}
over that for the classical \textsc{Svd} is of the order ${\mathcal
O}(l/n )$. Hence, the randomized algorithms can be essentially
faster than the traditional counterpart.

\subsection{Error estimates}
We will analyze the accuracy of the Algorithm \textsc{Rttls} in this part.
Before the main results, we introduce an important estimate in \cite{Halko}.
\begin{lemma} \emph{\textrm{\cite[Corollary 10.9]{Halko}}}
\label{lem3}
Suppose that $A \in \mathbb{R}^{m \times n}$ has
singular values $\sigma_1 \geq \sigma_2 \geq  \cdots.$ Choose a target rank
$k \geq 2$ and an oversampling parameter $p \geq 4,$ where $k+p \leq \min\{m, n\}.$
Draw an $n \times (k+p)$ standard Gaussian matrix $\Omega$ and let $Q$ be an
orthonormal matrix whose columns form a basis for the range of the sampled matrix
$A \Omega.$ Then
     \begin{eqnarray*}
                      \left\| A - QQ^\mathrm T A \right\|_2 & \leq & \left( 1 + 9 \sqrt{k+p} \sqrt{\min\{m, n\}}\right) \sigma_{k+1},
      \end{eqnarray*}
with failure probability at most $3p^{-p}.$
 \end{lemma}
From the process of Algorithm \textsc{Rttls}, we see that $U\Sigma
V^\mathrm T = QW\Sigma V^\mathrm T = QQ^\mathrm T C$ where we denote
$U = QW$ and $C = [A, ~b].$ So $\left\| C - U\Sigma V^\mathrm T
\right\|_2 = \left\| C - QQ^\mathrm T C \right\|_2.$ Hence we obtain
a good \textsc{Svd} approximation for $C$ with high probability.

Before studying the accuracy of the stochastic procedures in the algorithm,
we review the perturbation results given by Wei \cite[Theorem 4.1]{Wei}, which is stated in the
following slightly modified lemma.
\begin{lemma} \label{lem4}
Consider the \textsc{Tls} problem \eqref{eq1.1}. Let the
\textsc{Svd} for $A$ and $[A, ~b]$ be given as in the preliminaries.
Assume that for some $q \leq n,$ $\widetilde{\sigma}_q >
\sigma_{q+1}.$ Partition $V$ as in \eqref{Vpartion}, let
$\widehat{A} \in \mathbb{R}^{m \times n},~ \widehat{b} \in
\mathbb{R}^m,$ and $\left[\widehat{A},~\widehat{b}\right] =
\left[A,~b\right] + E$ with $\|E\|_2 \leq \frac16
\left(\widetilde{\sigma}_q - \sigma_{q+1}\right),$ and the
\textsc{Svd} for $\left[ \widehat{A}, ~~ \widehat{b} \right]$ be
$$ \widehat{U}^\mathrm T~\left[ \widehat{A}, ~~ \widehat{b} \right]~ \widehat{V}^\mathrm T = \widehat{\Sigma}.$$
Partition $\widehat{V}$ conformally with $V$ and replace $V_{ij}$ by $\widehat{V}_{ij}$
for $i,~j = 1, 2.$ Define $\widehat{x}_\textsc{ttls} = \left( \widehat{V}_{11}^\mathrm T \right)^\dag \widehat{V}_{21}^\mathrm T$
and $x_\textsc{ttls} = \left( V_{11}^\mathrm T \right)^\dag V_{21}^\mathrm T.$ When $x_\textsc{ttls} \neq {\bf 0},$ the following estimate holds:
     \begin{eqnarray*}
     \label{TTLS_perturb}
            \frac{\left\| x_\textsc{ttls} - \widehat{x}_\textsc{ttls} \right\|_2}{\left\| x_\textsc{ttls} \right\|_2}
             & \leq & \frac{12 \left(\|E\|_2 + \sigma_{q+1} \right)}{\widetilde{\sigma}_q - \sigma_{q+1}} \frac{\sigma_1}{\|b\|_2 - \sigma_{q+1}}.
      \end{eqnarray*}
 \end{lemma}
Using Lemma \ref{lem3} and Lemma \ref{lem4}, we get the estimate below.
\begin{theorem}
Assume $m \geq n+1.$
Assume that $[A, ~b]$ has singular values $\sigma_1, \dots, \sigma_{n+1}$ and
$A$ has singular values $\widetilde{\sigma}_1 \geq \cdots \geq \widetilde{\sigma}_n.$
Moreover, assume $\widetilde{\sigma}_q > \sigma_{q+1}$ with $q \leq n$  and
let $k$ be the target rank of $[A,~b]$ and let $x_\textsc{rttls}$ be the
approximate \textsc{Ttls} solution by performing
Algorithm \textsc{Rttls} with the Gaussian random matrix $\Omega \in \mathbb{R}^{n \times (k+p)}.$
If
     \begin{eqnarray*}
            \sigma_{k+1} \leq \frac{\widetilde{\sigma}_q - \sigma_{q+1}}{6 + 54 \sqrt{(k+p)n}},
      \end{eqnarray*}
then we have
     \begin{eqnarray}\label{error_RTTLS}
            \frac{\left\| x_\textsc{ttls} - x_\textsc{rttls} \right\|_2}{\left\| x_\textsc{ttls}\right\|_2} & \leq &
            \frac{12 \sigma_1 \left[ \left(1 + 9 \sqrt{(k+p)n}  \right)\sigma_{k+1} + \sigma_{q+1} \right]}{\left( \widetilde{\sigma}_q - \sigma_{q+1} \right) \left(\|b\|_2 - \sigma_{q+1} \right)}
      \end{eqnarray}
with failure probability at most $3 p^{-p}.$
More specifically, if $q = k$ is the numerical rank of $[A,~b],$ we get the bound below with probability not less than $1 - 3p^{-p}$
     \begin{eqnarray}\label{error_RTTLS_specific}
            \frac{\left\| x_\textsc{ttls} - x_\textsc{rttls} \right\|_2}{\left\| x_\textsc{ttls}\right\|_2} & \leq &
            \frac{12 \sigma_1 \left(2 + 9\sqrt{(k+p)n} \right)}{\widetilde{\sigma}_k \|b\|_2} \sigma_{k+1} + {\mathcal O}(\sigma_{k+1}^2).
      \end{eqnarray}
 \end{theorem}
\proof
Denote that $C = \left[A,~b\right]$ and $\widehat{C} = QQ^\mathrm T C.$ From Lemma \ref{lem3} and the assumption we know that
$$\left\| C - \widehat{C} \right\|_2 = \left\| C - QQ^\mathrm T C \right\|_2 \leq  \left(1 + 9 \sqrt{(k+p)n}  \right)\sigma_{k+1} \leq \frac16 \left(\widetilde{\sigma}_q - \sigma_{q+1}\right)$$
with probability not less than $1 - 3p^{-p}.$ Then applying Lemma \ref{lem4} we obtain \eqref{error_RTTLS}. For the specific cases,
if the numerical rank of $[A, ~b]$ is
$k = q,$ it means that $\sigma_{k+1}$ is very close to zero. The bound in \eqref{error_RTTLS} can be simplified to
     \begin{eqnarray*}
            \frac{\left\| x_\textsc{ttls} - x_\textsc{rttls} \right\|_2}{\left\| x_\textsc{ttls}\right\|_2} & \leq &
            \frac{12 \sigma_1 \left(2 + 9 \sqrt{(k+p)n}  \right) }{\left( \widetilde{\sigma}_k - \sigma_{k+1} \right) \left(\|b\|_2 - \sigma_{k+1} \right)}
            \sigma_{k+1}.
      \end{eqnarray*}
 Consider the Taylor expansion for the function
 $f(x) = 1 / \left[ \left( \widetilde{\sigma}_k - x \right) \left(\|b\|_2 - x \right) \right]$ at $x = 0,$ we obtain that
 $$f(\sigma_{k+1}) = \frac{1}{ \widetilde{\sigma}_k  \|b\|_2 } + \left( \frac{1}{\widetilde{\sigma}_k ^2 \|b\|_2} + \frac{1}{\widetilde{\sigma}_k \|b\|_2 ^2}\right)\sigma_{k+1} + {\mathcal O}(\sigma_{k+1}^2).$$
 Substituting this equation into the above inequality directly, we can get \eqref{error_RTTLS_specific}.
\proofend

We point out that the assumption for $\sigma_{k+1}$ usually holds for the ill-conditioned cases, where $k = q$ is the numerical rank and we treat the other smaller singular values as zeros. During these cases, the upper bound \eqref{error_RTTLS_specific} is of order ${\mathcal O}(\sigma_{k+1})$ and hence the relative error of the solution from \textsc{Rttls} and the solution from \textsc{Ttls} is small.

\section{Numerical examples}
In this section we give numerical examples to verify the perturbation
bounds and our randomized total
least squares algorithms (\textsc{Rtls} and \textsc{Rttls}). The
following numerical tests are performed via \textsc{Matlab} R2010a
in a laptop with Intel Core i5 by using double precision.

\subsection{Perturbation bounds}
We compare our upper bounds (\ref{err1}) and (\ref{Corollary}) with those derived in
\cite{Baboulin,Jia,Zhou}. We will see that these three are equal,
and ours are sharper.

{\bf Example I.}
 In this example \cite[Example 1]{Baboulin}
we consider the \textsc{Tls} problem $Ax\approx b$, where $[A, ~b]$
is defined by
$$[A, ~b]=Y\left[\begin{array}{c}
                         D \\
                         {\bf 0}
                          \end{array} \right]Z^\mathrm T \in \mathbb{R}^{m\times (n+1)}, Y=I_m -2yy^\mathrm T, Z=I_{n+1}-2zz^\mathrm T,$$
where $y\in \mathbb{R}^m$ and $z \in \mathbb{R}^{n+1}$ are random unit vectors, $D={\rm diag}(n,n-1,\cdots,1,1-\epsilon_p)$ for a given parameter $\epsilon_p$. The quantity $\widetilde{\sigma}_n-\sigma_{n+1}$  measures the distance of our problem to nongenericity and, due to the interlacing property, we have in exact arithmetic $$\widetilde{\sigma}_n - \sigma_{n+1} \leq \sigma_n - \sigma_{n+1}=\epsilon_p.$$
We consider a random perturbation $\|[\delta A,~\delta b]\|_F = 10^{-10}.$
We take $m=100,~n=40$ in this example and denote $\Delta = \frac{\|[\delta A, ~\delta b]\|_F}{\|[A, ~b]\|_F}$.

\begin{table}\centering

\begin{tabular}{ccccccc}\hline
 $\epsilon_p$ & $\frac{\|\hat{x}-x\|_2}{\|x\|_2}$ & $K_\textsc{zlwq} \Delta$ & $K_\textsc{bg} \Delta$ & $K_\textsc{lj} \Delta$ & (\ref{err1}) & (\ref{Corollary})\\
\hline
 9.99976032E-1 &2.6233E-11  &1.5815E-09 &1.5815E-09 &1.5815E-09  &2.8145E-10 &3.9565E-10\\
 9.99952397E-5 &3.8714E-07   &1.1343E-05 &1.1343E-05 &1.1343E-05  &3.2472E-06 &3.8752E-06\\
\hline
\end{tabular}
\caption{\label{tab:boundsExampleI}Comparisons of forward error and
upper bounds for a perturbed \textsc{Tls} problem.}
\end{table}

In Table \ref{tab:boundsExampleI}, we compare the exact relative
error with the upper bounds (\ref{err1}) and the above bounds
derived in \cite{Baboulin,Jia,Zhou}.
Without considering the computational cost, we can see that
the numerical results of the three condition numbers in \cite{Baboulin,Jia,Zhou} are the same.
We observe that our bounds are  sharp and  smaller than the bounds derived in the literature.

{\bf Example II.} Consider the second example from \cite[p.
42]{VanHuffel}, where
$$A = \left[\begin{array}{cccc}
                         m-1 & -1 & \cdots & -1 \\
                         -1 & m-1 & \cdots & -1 \\
                         \vdots & & &           \\
                         -1 & -1 & \cdots & m-1 \\
                         -1 & -1 & \cdots & -1 \\
                         -1 & -1 & \cdots & -1
                          \end{array} \right] \in \mathbb{R}^{m \times (m-2)}, \quad b =
                          \left[\begin{array}{c}
                         -1 \\
                         -1 \\
                         \vdots\\
                         -1 \\
                         m-1  \\
                         -1
                          \end{array} \right] \in \mathbb{R}^m.$$
The exact solution of the \textsc{Tls} problem $Ax \approx b$ is $x
= - [1,1,\ldots,1]^\mathrm T$ and $\sigma_{n+1}=\sqrt{m}, ~
\widetilde{\sigma}_n = \sqrt{2m}.$ We consider the same random
perturbations as the Example I. The results are listed in Table
\ref{tab:boundsExampleII}.
\begin{table}\centering

\begin{tabular}{cccccccc}\hline
 $m$ & $\frac{\|\hat{x}-x\|_2}{\|x\|_2}$ & $K_\textsc{zlwq} \Delta$
 & $K_\textsc{bg} \Delta$ & $K_\textsc{lj}\Delta$ & (\ref{err1}) & (\ref{Corollary})\\
\hline
 100 &1.1553E-13  &1.0152E-11 &1.0152E-11 &1.0152E-11 &4.7548E-12 &4.6281E-12  \\
 250 &2.5302E-14  &6.3627E-12  &6.3627E-12 &6.3627E-12 &1.9277E-12  &1.9001E-12  \\
\hline
\end{tabular}
\caption{\label{tab:boundsExampleII}Comparisons of forward error and
upper bounds for a perturbed \textsc{Tls} problem. }
\end{table}
From the experience of computing, we also find that the bound in
\cite{Zhou} is quite impractical for computing, since
\textsc{Matlab}  will be out of memory on our Microsoft Windows
operating system.

\subsection{Numerical experiments for randomized algorithms}
In this subsection, we apply Algorithm \textsc{Rtls}, Algorithm
\textsc{Rttls} and Algorithm \textsc{Arttls} to Example I, Example II
and some cases in Hansen's Regularizaton Tool \cite{Hansen94}.
We will compare the
computational time and solution accuracy of our new randomized
\textsc{Tls} algorithms with the traditional algorithms.

\subsubsection{Algorithm \textsc{Rtls} on well-conditioned cases}

\begin{table}\centering

\begin{tabular}{llccccl}   \hline
 &   Matrix size   & $\mathrm{Cond}(A)$  & $\mathrm{Cond}([A, ~b])$ &  $\mathrm{Time}_\textsc{Tls}$ & $\mathrm{Time}_{\textsc{Rtls}}$ & $\mathrm{Err}_\textsc{Rtls}$
 \\ \hline 
 Example I  &  $m=500$   & 2.00E+2 & 8.34E+6 & 0.0698  & 0.0110 &  6.48E-10  \\
            &  $m=1000$  & 4.00E+2& 1.67E+7 & 0.5643    & 0.0941  & 1.06E-10  \\
            &  $m=5000$  & 2.00E+3& 8.34E+7 & 35.021    & 2.7903  & 2.40E-09  \\
 Example II  &  $m=500$         & 15.8 & 22.4   &  0.2063   & 0.0402 & 5.53E-02  \\
            &  $m=1000$         & 22.4 & 31.6  & 1.3648 &  0.2706 & 4.09E-02  \\
            &  $m=5000$         & 50.0 & 70.7  & 154.51 &  23.345 & 1.88E-02  \\
 Deriv2     &  $m=500$          & 3.04E+5 & 3.33E+5 & 0.3180 &   0.0477  &  5.34E-05 \\
            &  $m=1000$         &1.22E+6 & 1.33E+6 &  1.7050    &  0.4033  & 6.56E-04   \\
            &  $m=5000$         &3.04E+7 & 3.33E+7 &  157.64    &  39.114  & 5.15E-01   \\ \hline
\end{tabular}
\caption{\label{tab:RTLS} Tests on  \textsc{Rtls}  ($l=10$).
 $\mathrm{Time}_\textsc{Tls}$
and $\mathrm{Time}_{\textsc{Rtls}}$ are the computational times (in
seconds) for the algorithms \textsc{Tls} and \textsc{Rtls}
respectively. The relative error $\mathrm{Err}_\textsc{Rtls} = ||
x_\textsc{tls} - x_\textsc{rtls} ||_\infty /
\|x_\textsc{tls}\|_\infty.$
 }
\end{table}
For the case Example I in Table \ref{tab:RTLS}, we choose
$\epsilon_p$ =9.99976031e-1, and set $n = \frac{2}{5} m$. The
solution $x_\textsc{tls}$ is computed by \eqref{x}, while
$x_\textsc{rtls}$ is obtained by Algorithm \textsc{Rtls}. Denote the
relative error $\mathrm{Err}_\textsc{Rtls} = || x_\textsc{tls} -
x_\textsc{rtls} ||_\infty / \|x_\textsc{tls}\|_\infty.$
 The corresponding execution time $\mathrm{Time}_\textsc{Tls}$
and
$\mathrm{Time}_{\textsc{Rtls}}$ are measured by the \textsc{Matlab}
tic-toc pairs in seconds. From Table \ref{tab:RTLS}, we can see that
our \textsc{Rtls} algorithm on large matrices outperforms the
traditional counterpart according to computational time, while the
accuracy of solutions of two methods is comparable. For the small
matrices, the advantage of \textsc{Rtls} will not be so obvious. For
an ill-conditioned matrix, \textsc{Matlab} reports inaccuracy
warning due to the ill-conditioned linear system in Step 2 and 4 of
Algorithm \textsc{Rtls}. Even for the ill-conditioned case
\textsc{Deriv2}, Algorithm \textsc{Rtls} can still give approximate
solution with good accuracy. But for the very ill-conditioned cases,
we need Algorithm \textsc{Rttls}.



\subsubsection{Examples based on \textsc{Tls}-Prony modeling}
\begin{figure}[htbp] \centering
\includegraphics[width=\textwidth]{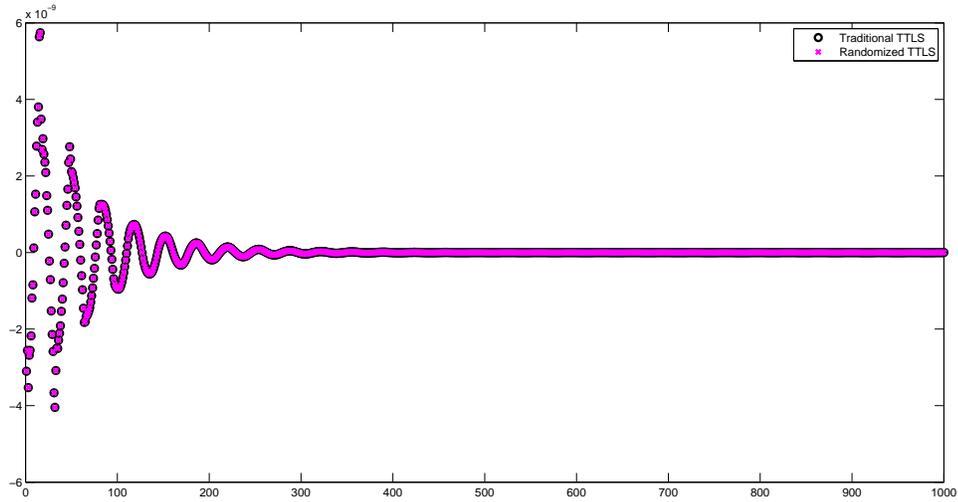}
\caption{\label{Fig:Prony_1000}Computed solutions for the
true \textsc{Tls} solution and the randomized one. }
\end{figure}

The \textsc{Tls} approach is a promising method in the field of signal
processing. Rahman and Yu \cite{Rahman} presented a method for frequency estimation
using \textsc{Tls} for solving the linear prediction equation.
The problem here is taken from \cite{Majda}. We first consider a set
of linear prediction equations. Assume $a_j = \left[y_{j-1},\ldots, y_{j+m-2} \right]^\mathrm T$
where $y_l = \sum_{j=1}^p c_j z_j^l,~z_j = \exp\left(\lambda_j T\right),~ j=1,\ldots,p.$ The $\lambda_j$'s
and $c_j$'s are to be determined. Furthermore, assume $c_j$ and $z_j$ are nonzeros and $z_j$'s are distinct
for $j = 1,\ldots, p.$ Let $A_n = \left[a_1, \ldots, a_n\right],~ b_n = -a_{n+1}$ and consider the linear system
     \begin{eqnarray} \label{Prony}
                       A_n x & = b_n.
      \end{eqnarray}
Assume $m \geq n,~ m \geq p.$ It is known \cite{weiLAA} that $\mbox{rank}\left(A_n\right) = \min \{n,p\}.$ So if $n \geq p,$ then \eqref{Prony}
is compatible. For any solution $x = \left(\alpha_0, \alpha_1,\ldots, \alpha_{n-1} \right)^\mathrm T,$ construct a polynomial
$$P_n(z) = z^n + \alpha_{n-1} z^{n-1} + \cdots + \alpha_1 z + \alpha_0,$$
then we know that $P_n$ has zeros $z_1, \ldots, z_p.$
We choose $\lambda_j$ and $c_j$ as in Table \ref{tab:prony}.
\begin{table}\centering
\begin{tabular}{cc}\hline
  $\lambda_j$  & $c_j$
 \\ \hline
  $-0.082 \pm 0.926 i$ &  1  \\
  $-0.147 \pm 2.874 i$ &  1  \\
  $-0.188 \pm 4.835i$ &  1 \\
  $-0.220 \pm 6.800i$ &  1 \\
  $-0.247 \pm 8.767i$  &  1 \\
  $-0.270 \pm 10.733i$ &  1 \\
  \hline
\end{tabular}
\caption{\label{tab:prony}Six pairs of poles and residues. }
\end{table}
In this example $T = 0.2,~ m = 2000, ~ p = 12, ~n = 1000$ are used
and we compare the \textsc{Ttls} with the \textsc{Rttls} where the
sampling size is chosen as $l = p+1$. The plots for the solutions
are shown in Figure \ref{Fig:Prony_1000} and the infinity norm
relative error is $6.7623e-8$, while the time for \textsc{Tls} using
partial \textsc{Svd} and \textsc{Rttls} are 0.8924 seconds and
0.0333 seconds respectively.

\subsubsection{Algorithm \textsc{Rttls} on ill-conditioned cases}

Our ill-conditioned cases are taken from Hansen's Regularizaton
Tools \cite{Hansen94}. For example, the case \textsc{Shaw} is
generated by the command $\left[ \bar{A}, ~\bar{b}, ~x_\mathrm{true} \right]=
\textsc{Shaw} (m)$. Then noises are added to $\bar{A}$ and
$\bar{b}$. Suppose that $\delta$ is the relative noise level. We
define
\begin{eqnarray*}
b =  \bar{b} + \delta \left\| \bar{b} \right\|_2 \frac{\zeta}{\left\| \zeta \right\|_2}, \qquad
A =  \bar{A} + \delta \left\| \bar{A} \right\|_F \frac{Z}{\| Z \|_F},
\end{eqnarray*}
where $\zeta$ is a random vector, $\zeta= 2*\mathrm{rand}(m,1)-1$;
$Z$ is a random matrix, $Z = 2*\mathrm{rand}(m)-1$. It is easy to
verify that
$$ \frac{\left\| b - \bar{b}\right\|_2}{\left\| \bar{b} \right\|_2} = \frac{\left\| A-\bar{A} \right\|_F}{\left\| \bar{A} \right\|_F} = \delta.$$
Then we seek the total least squares solution of $A x \approx b$.

\begin{figure}[htbp] \centering
\includegraphics[width=\textwidth]{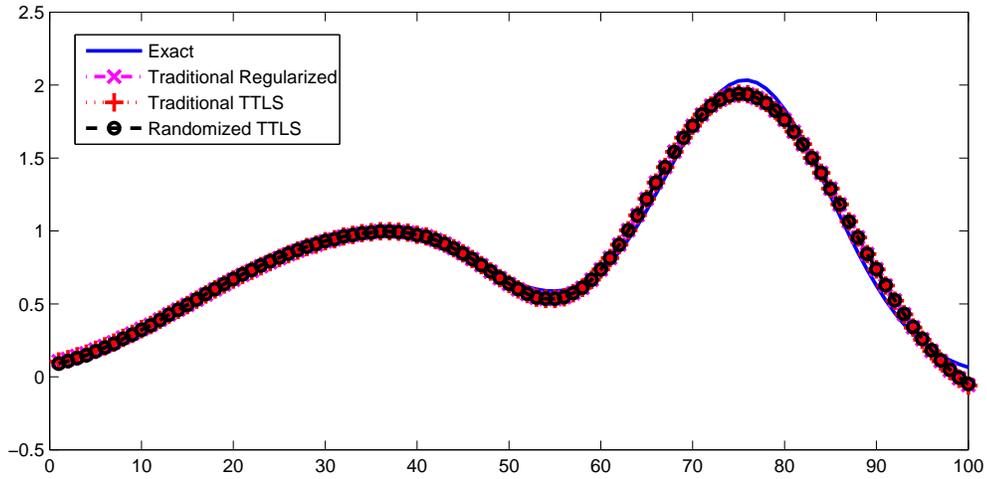}
\caption{\label{Fig:shaw_n100_noise_3}Computed solutions for the
case \textsc{Shaw} of size $m$=100 with relative noise level
$\delta$=1E-3. }
\end{figure}

\begin{table}\centering

\begin{tabular}{cccccl}\hline
  $\delta$  & $k$ &  $\mathrm{Time}_\textsc{Ttls}$ & $\mathrm{Time}_\textsc{Ptls}$ & $\mathrm{Time}_\textsc{Rttls}$ & $\mathrm{Err}_\textsc{Rttls}$
 \\ \hline
  1E-1 &  3 &   0.0123 & 0.2230 &  0.0039  & 8.04E-3  \\
  1E-2 &  5 &   0.0133 & 0.2365 &  0.0039 &  8.92E-4 \\
  1E-3 &  7 &   0.0119 & 0.2455 &  0.0075  & 1.59E-3  \\
  1E-4 &  8 &   0.0114 & 0.2424 &  0.0039 &  3.76E-4 \\
            \hline
\end{tabular}
\caption{\label{tab:rrtls_shaw}Tests on  \textsc{Shaw} with
different relative noise levels. The relative error
$\mathrm{Err}_\textsc{Rttls} = || x_\textsc{rttls} - x_\textsc{ttls}
||_\infty / \|x_\textsc{ttls}\|_\infty$. Algorithm \textsc{Rttls} is
substantially faster than  \textsc{Ttls} and \textsc{Ptls}.
}
\end{table}
We first test Algorithm \textsc{Rttls} on the $100 \times 100$
matrix \textsc{Shaw} with different relative noise levels $\delta$.
The results are given in Table \ref{tab:rrtls_shaw}. The truncation
parameter $k$ is estimated by the randomized algorithm with
\textsc{Gcv} and \textsc{Tsvd} \cite{Xiang}. After the determination
of parameter $k$, the computational time for implementing Algorithm
\textsc{Ttls} and Algorithm \textsc{Rttls} is recorded in
$\mathrm{Time}_\textsc{Ttls}$ and $\mathrm{Time}_\textsc{Rttls}$
respectively. And $\mathrm{Time}_\textsc{Ptls}$ denotes the time
cost in \textsc{Ttls} using Lanczos bi-diagonalization based partial
\textsc{Svd}. Here we denote the relative error
$\mathrm{Err}_\textsc{Rttls} = || x_\textsc{ttls} - x_\textsc{rttls}
||_\infty / \|x_\textsc{ttls}\|_\infty.$ From our computing, we see
that the computed solutions of Algorithm \textsc{Ptls} are almost
the same as those of Algorithm \textsc{Ttls}, and hence the relative
errors for the solutions of \textsc{Ptls} which we denote as
$\mathrm{Err}_\textsc{Ptls} = \| x_\textsc{ptls} - x_\textsc{ttls}
\|_\infty / \|x_\textsc{ttls}\|_\infty$ are much smaller. 
Here we ignore the error $\mathrm{Err}_\textsc{Ptls}$ and do not
list it in the table.
 From this table
we can see that the results of Algorithm \textsc{Rttls} are very
close to those of traditional \textsc{Ttls} even for the relative
noise level as large as 10\%. According to the computational time,
Algorithm \textsc{Ptls} does not show obvious advantages over
\textsc{Ttls} for small size cases, while Algorithm \textsc{Rttls}
is substantially faster than the traditional \textsc{Ttls}. The
computed solutions for the case where the relative noise level
$\delta$=1E-3 are presented in Figure \ref{Fig:shaw_n100_noise_3}.

\begin{figure}[!htbp]\centering
\subfigure[][\textsc{Baart}]
   {\includegraphics[scale=0.18]{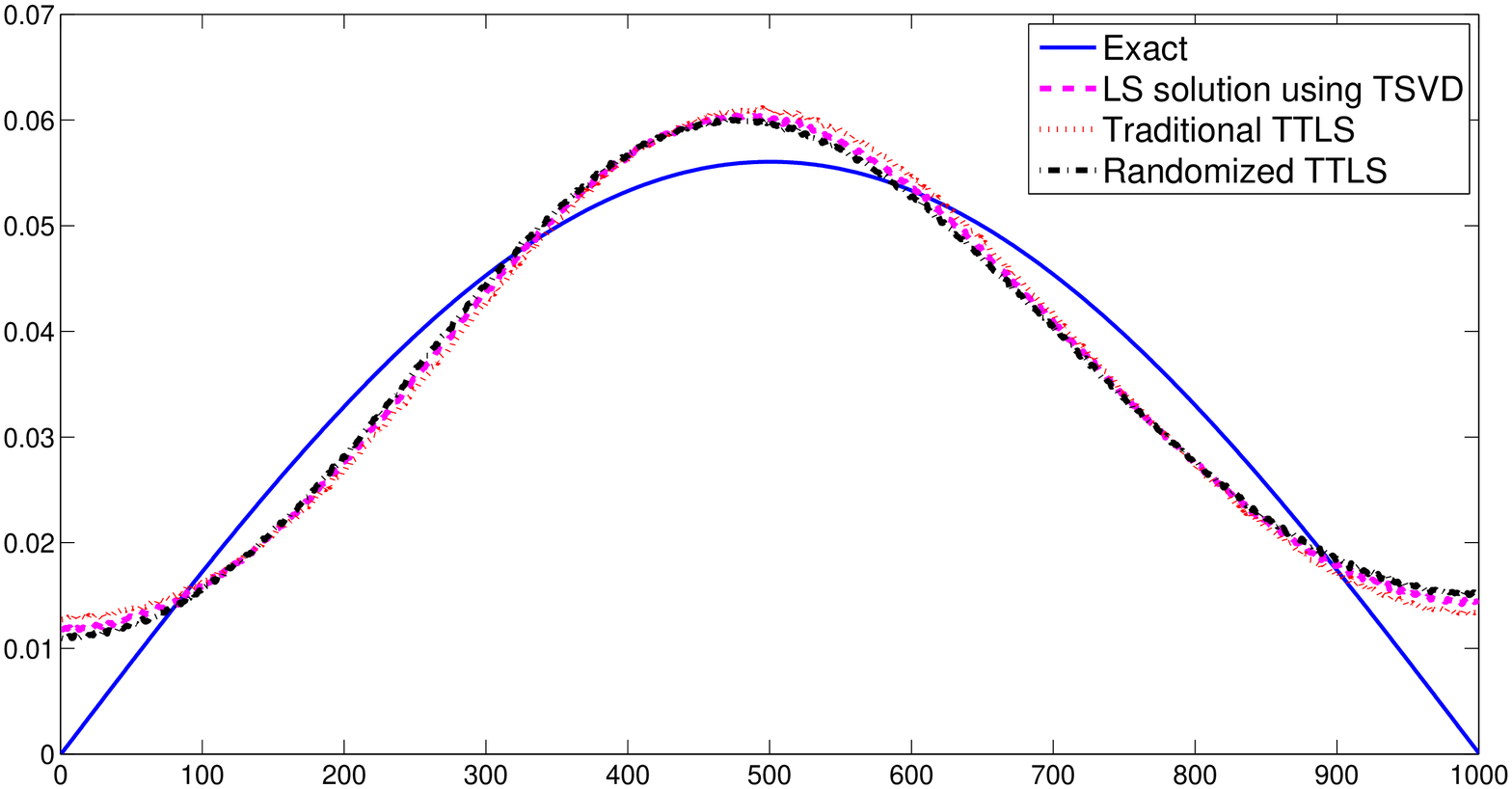}\label{Fig_Baart}}
\subfigure[][\textsc{Deriv2}]
   {\includegraphics[scale=0.18]{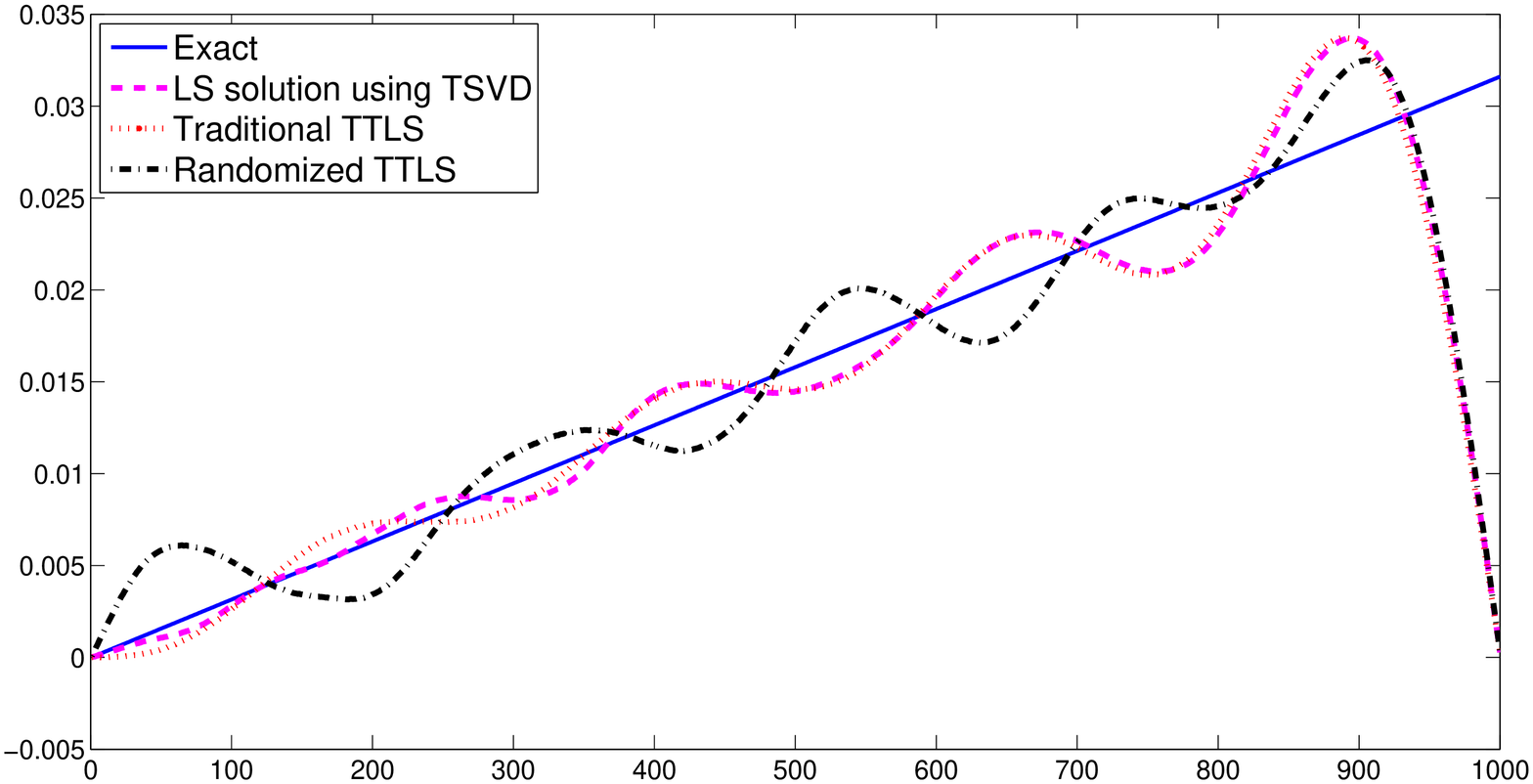}\label{Fig_Deriv2}}
\subfigure[][\textsc{Foxgood}]
   {\includegraphics[scale=0.18]{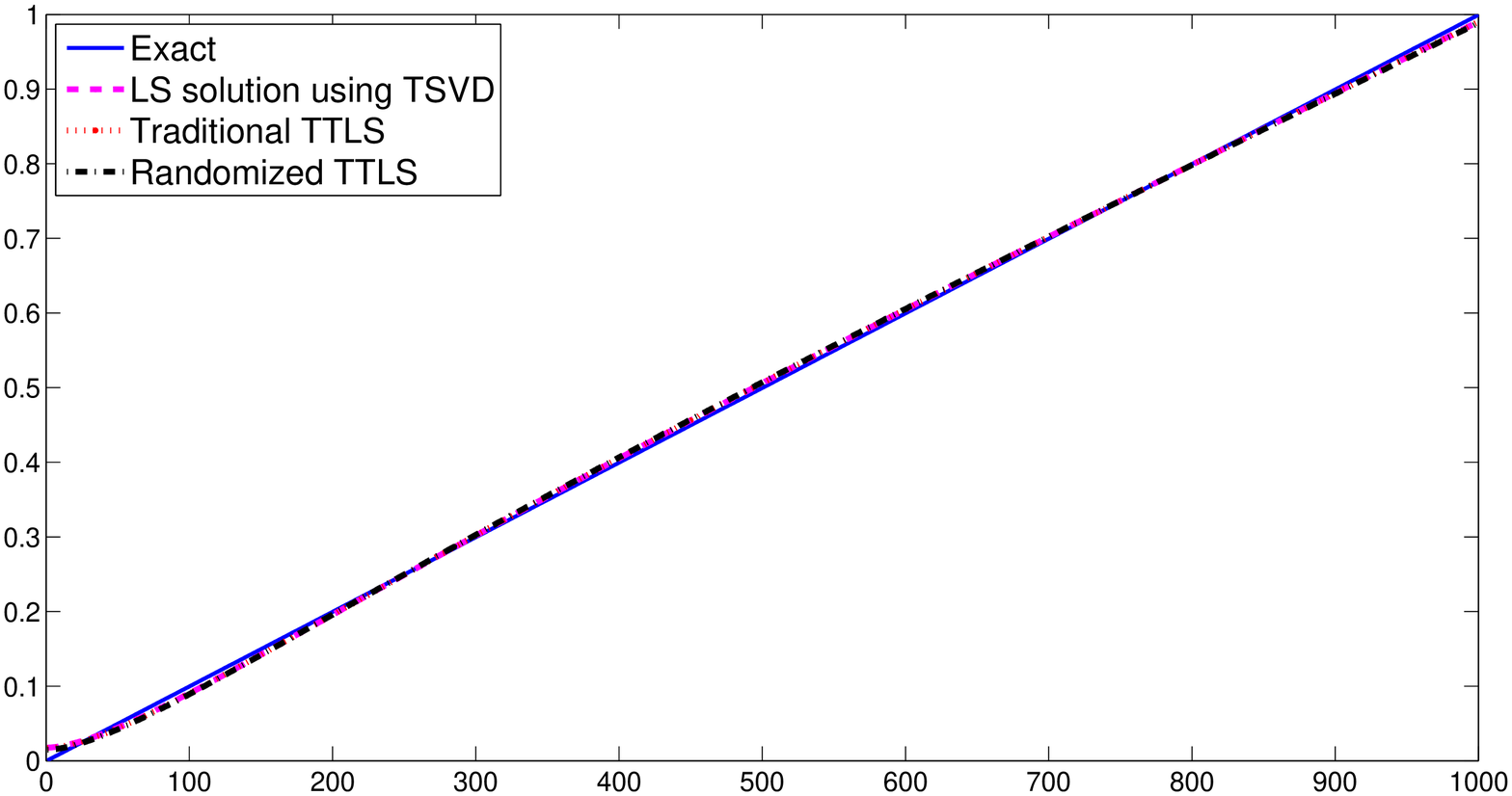}\label{Fig_Foxgood}}
\subfigure[][\textsc{Gravity}]
   {\includegraphics[scale=0.18]{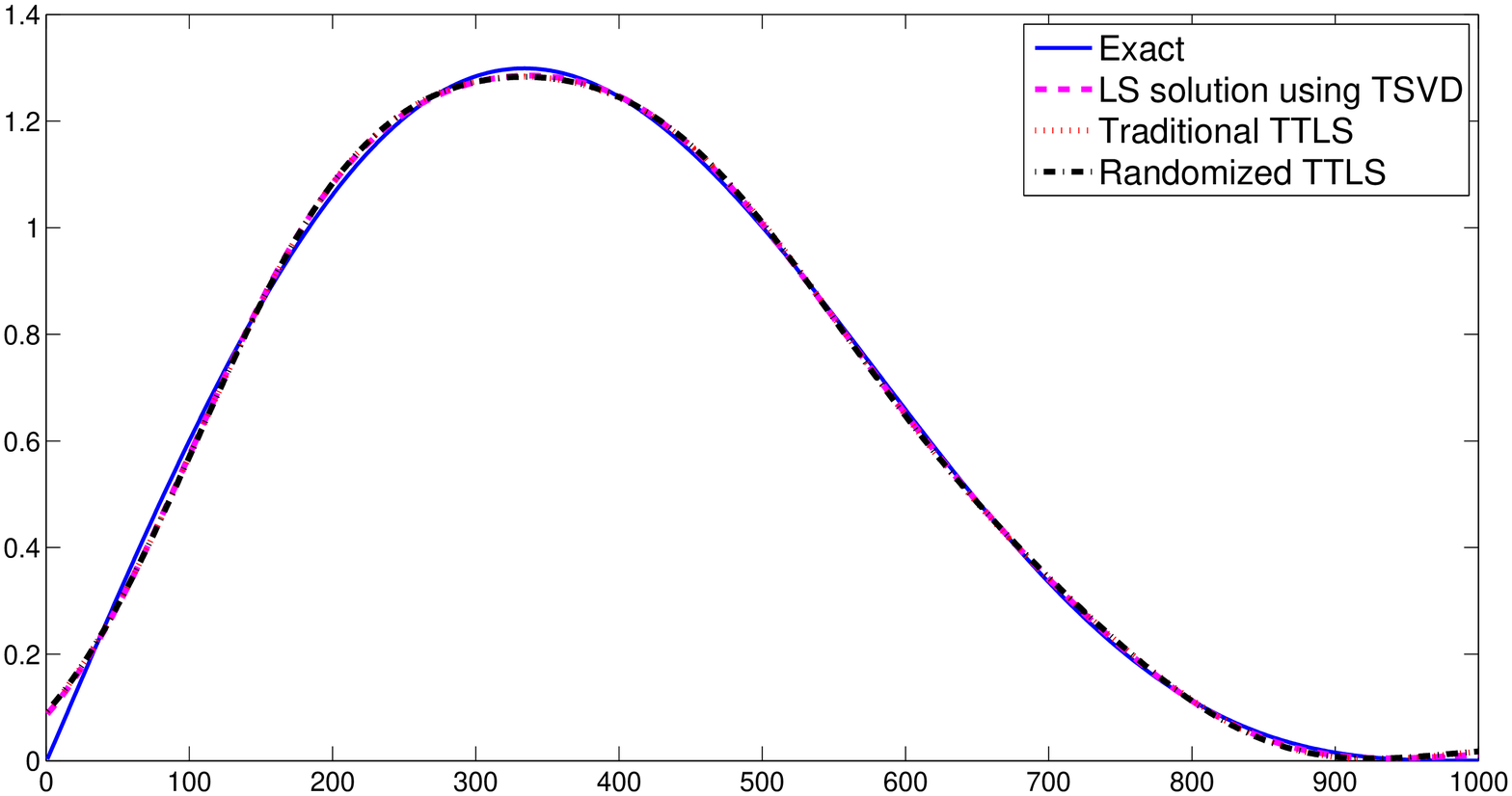}\label{Fig_Gravity}}
\subfigure[][\textsc{Heat}]
   {\includegraphics[scale=0.18]{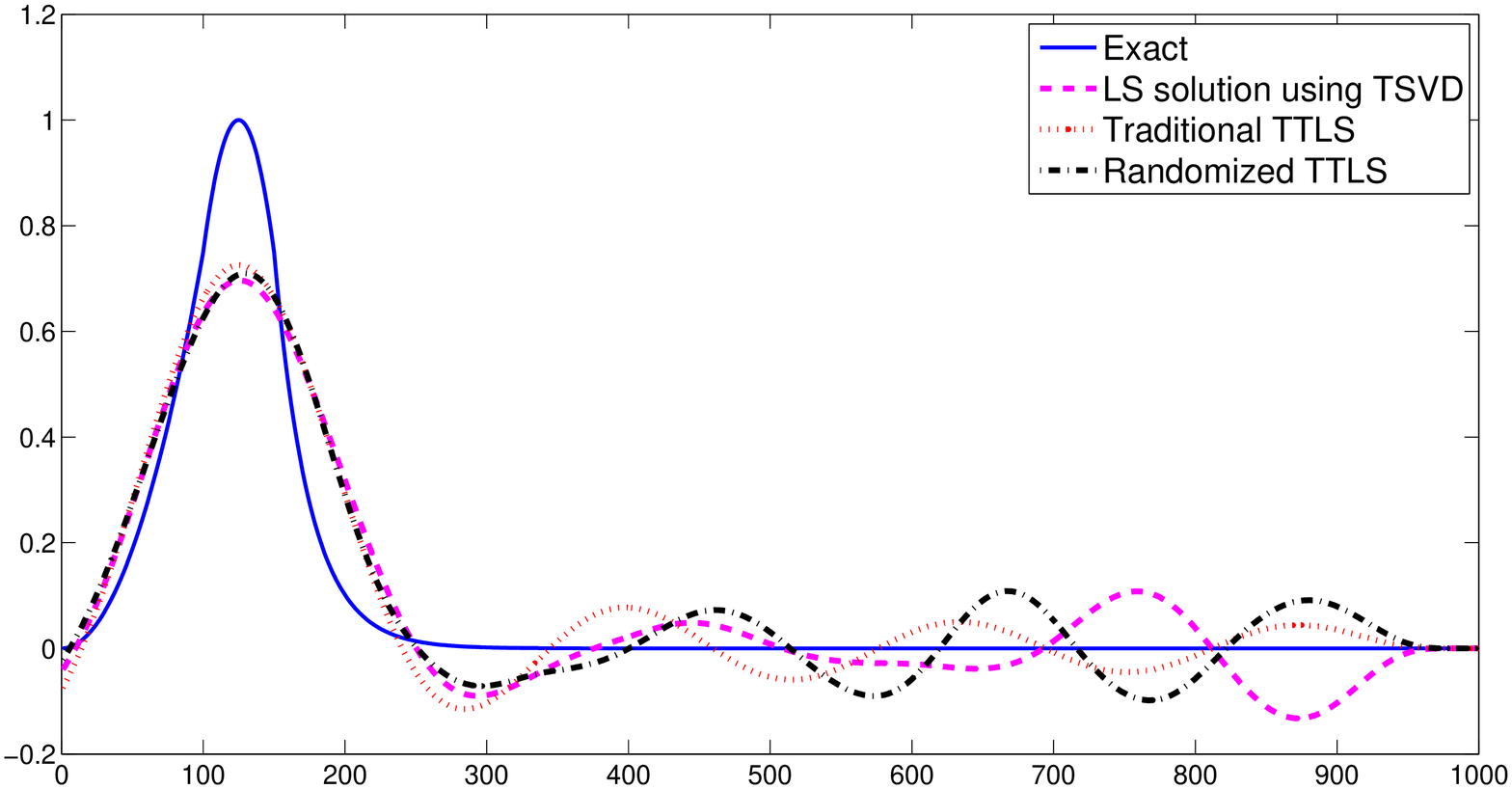}\label{Fig_Heat}}
\subfigure[][\textsc{i\_Laplace}]
   {\includegraphics[scale=0.18]{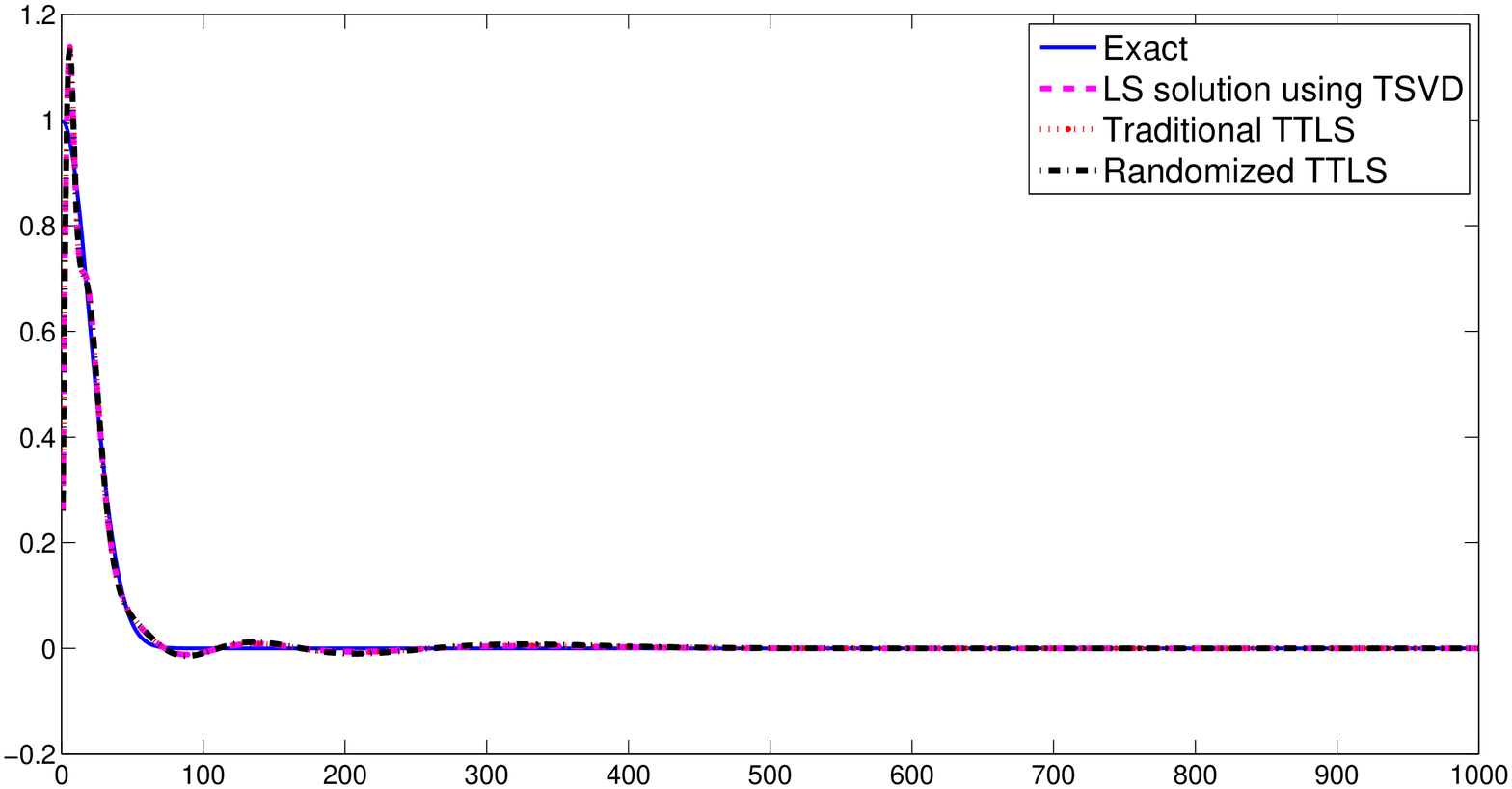}\label{Fig_ILaplace}}
\subfigure[][\textsc{Phillips}]
{\includegraphics[scale=0.18]{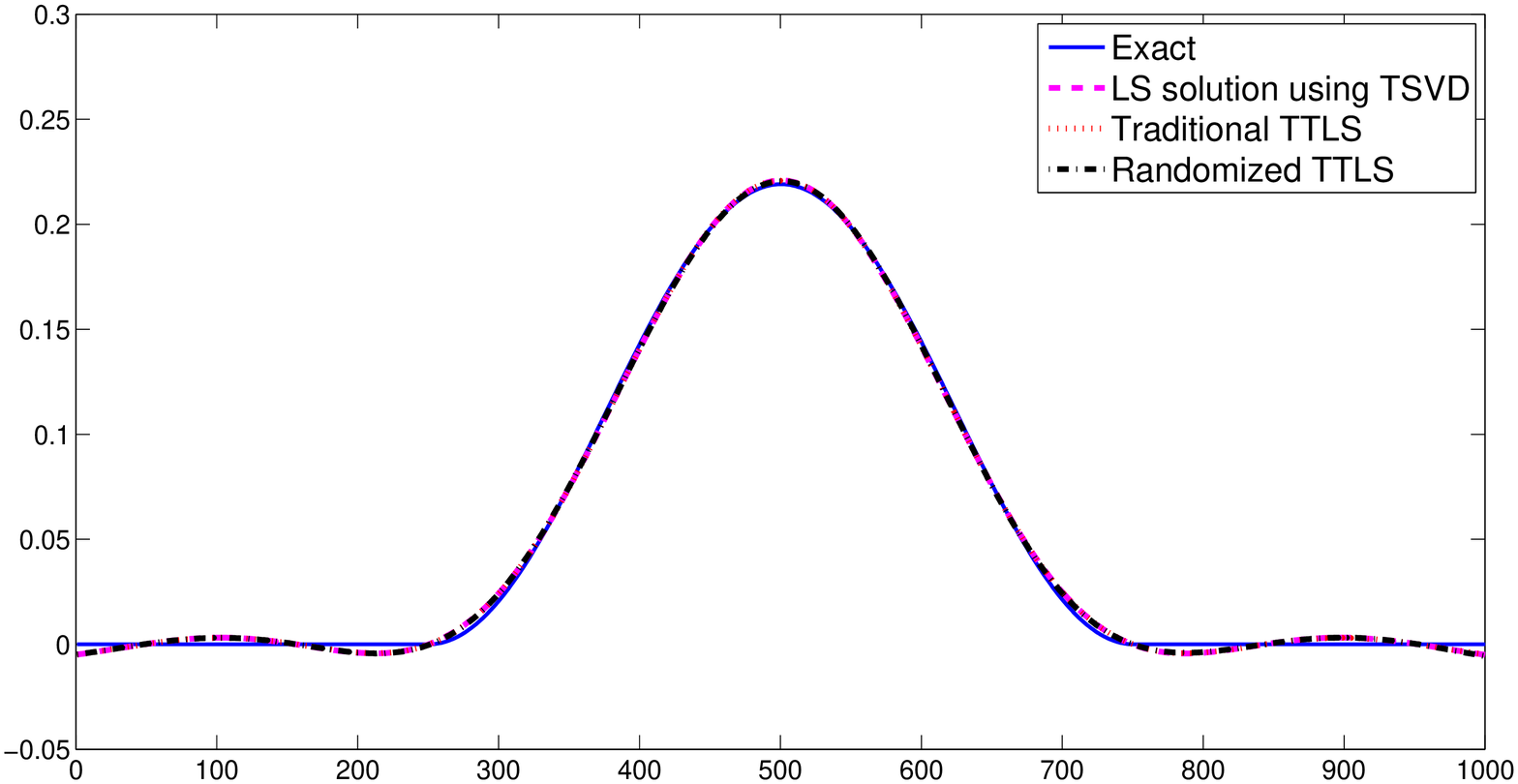}\label{Fig_Phillips}}
\subfigure[][\textsc{Shaw}]
   {\includegraphics[scale=0.18]{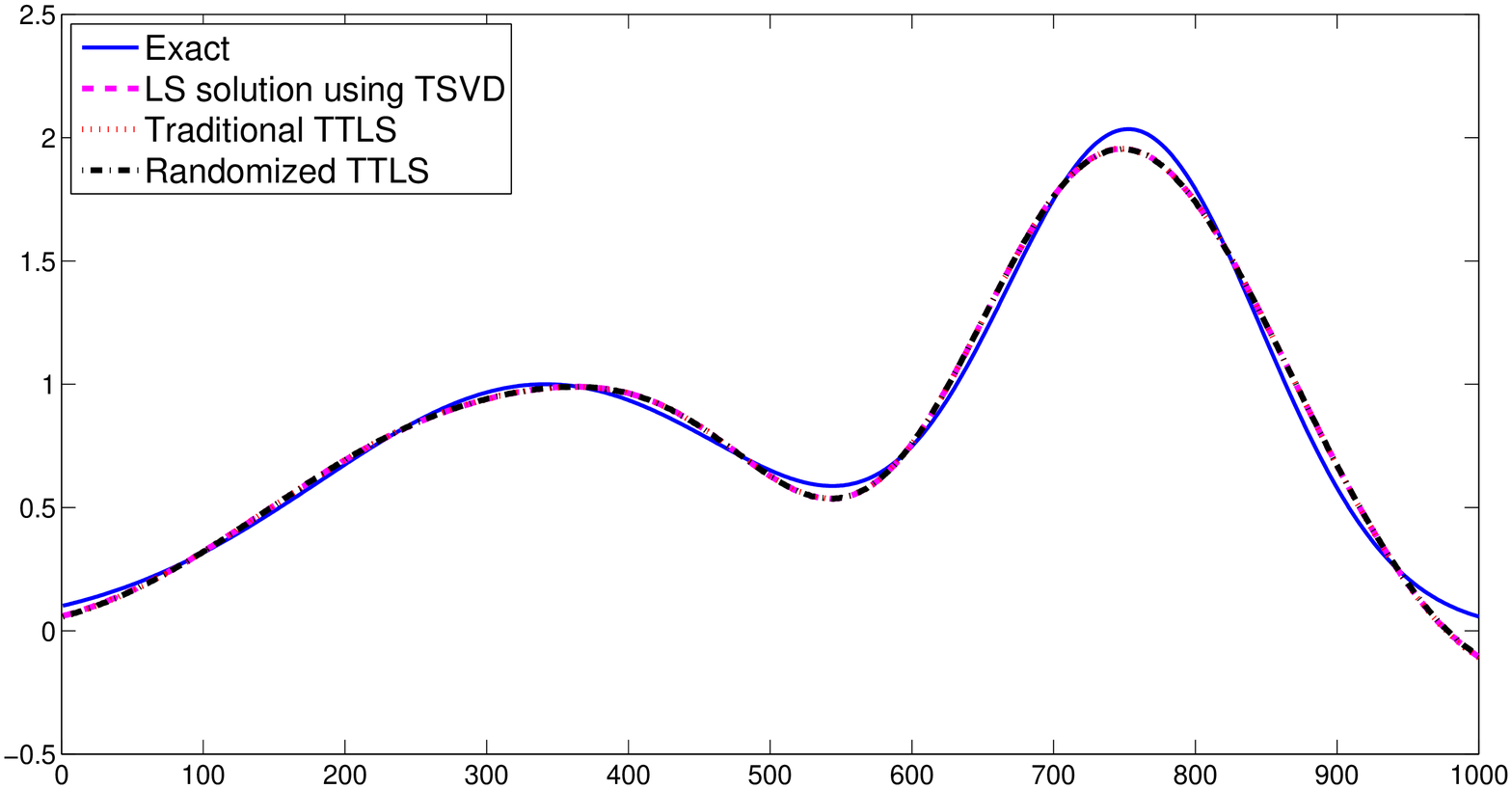}\label{Fig_Shaw}}

\caption{ \label{fig:8plots} \textsc{Rttls}  for ill-conditioned
cases of size $m=1000$ with relative noise level $\delta$=1E-3.  }
\end{figure}

We then test Algorithm \textsc{Rttls} on larger matrices. We set the
parameter for sampling size $l=10$ and the relative noise level
$\delta$=1E-3 for all cases. The results are given in Table
\ref{tab:rttls}. The marker $*$ in Table \ref{tab:rttls} means we
cannot load the example I\_Laplace on our computer when the size $n
= 5000.$ Obviously, \textsc{Ptls} can be much faster than the
\textsc{Ttls} when the size of the matrix becomes larger. But it is
still not as efficient as the randomized one because of its
\textsc{Blas}-2 operations. The randomized strategy can greatly
speed up the classical Algorithm \textsc{Ttls}.
  The advantage of our Algorithm \textsc{Rttls} is
more obvious when we test the larger matrices. The plots of the
computed solutions are given in Figure \ref{fig:8plots}.

\begin{table}\centering
\begin{tabular}{llccccc}\hline
 &   Matrix size     & $k$ &  $\mathrm{Time}_\textsc{Ttls}$ & $\mathrm{Time}_\textsc{Ptls}$ & $\mathrm{Time}_\textsc{Rttls}$ & $\mathrm{Err}_\textsc{Rttls}$
 \\ \hline
  Baart     &  $m=100$   &   4 & 0.0153  & 0.2907   & 0.0040 & 6.43E-3 \\
            &  $m=1000$   &  4 & 1.7561  &0.2664   & 0.0143 & 6.53E-3\\
            &  $m=5000$   &  4 & 176.47  &1.5042 & 0.2471 & 5.86E-3\\
  Deriv2     &  $m=100$   &  6  & 0.0130 &  0.2604  &0.0040 & 1.39E-2 \\
            &  $m=1000$   &  7  & 1.6727 &  0.3589  & 0.0148 & 6.96E-2 \\
            &  $m=5000$   &  9  & 170.40 &1.8398  & 0.2506 & 1.20E-2 \\
  Foxgood     &  $m=100$  & 2   &0.0129 & 0.2691  & 0.0037  & 4.60E-6  \\
            &  $m=1000$   & 3   &1.7638 &  0.2795  & 0.0143 & 5.09E-4 \\
            &  $m=5000$   & 3   &171.88 & 1.1383  & 0.2227 & 1.14E-4 \\
  Gravity    &  $m=100$   & 7   &0.0152  & 0.4679  & 0.0039 & 1.91E-3 \\
            &  $m=1000$   & 8  &1.7214   & 0.2963 & 0.0147 & 6.70E-3 \\
             &  $m=5000$   & 9  &183.92 &2.2156 & 0.3014 & 3.16E-2 \\
  Heat    &  $m=100$      & 8  &0.0107   & 0.2283   &  0.0041  & 7.33E-2 \\
            &  $m=1000$   & 9  &1.7172   & 0.3963   & 0.0163   & 3.93E-2 \\
            &  $m=5000$   & 9  &165.56   &1.4494   & 0.2551  & 8.15E-2 \\
  I\_Laplace    &  $m=100$  & 8 &0.0182   &0.2972   &0.0056     &2.22E-4   \\
            &  $m=1000$   & 9  &3.4499   &0.5609    &0.0410    &1.83E-2    \\
            &  $m=5000$   & *  &*   &*    &*    &*    \\
  Phillips    &  $m=100$  & 7 &0.0109  & 0.2476  & 0.0038 &  1.66E-3 \\
            &  $m=1000$   & 7  &2.3627 & 0.2804   &0.0137 & 2.24E-3 \\
            &  $m=5000$   & 7  &174.74 & 1.1844   &0.2194 & 6.08E-3 \\
           \hline
\end{tabular}
\caption{\label{tab:rttls} Algorithm \textsc{Rttls} on
ill-conditioned cases. 
}

\end{table}

\subsubsection{Test on Adaptive Algorithm \textsc{Arttls}}
We still use the examples from Hansen's Regularizaton Tools
\cite{Hansen94} and test the case with matrix size $m = n =1000.$
Here we set $r = 7$ in the Algorithm \textsc{Arttls}. Different
tolerances generate different $j$'s in Algorithm \textsc{Arttls}. So
we tried several $\epsilon$'s to make sure that $k$'s in Algorithm
\textsc{Rttls} and $j$'s in Algorithm \textsc{Arttls} are close, and
then the comparisons for the relative errors and time are
reasonable. The performance for Algorithm \textsc{Arttls} is shown
in Table \ref{tab:arttls}. In the table,
$\mathrm{Err}_\textsc{Arttls}$ denotes the relative error $||
x_\textsc{ttls} - x_\textsc{arttls} ||_\infty /
\|x_\textsc{ttls}\|_\infty$ and $\mathrm{Time}_\textsc{Arttls}$
represents the time cost for Algorithm \textsc{Arrtls}. It is clear
that Algorithm \textsc{Arttls} can still give good accuracy with
less computational time than the traditional one under the fixed
precision.
\begin{table}\centering

\begin{tabular}{llccccccc}\hline
 & $k$ &  $\epsilon$ & $j$ & $\mathrm{Time}_\textsc{Ptls}$ & $\mathrm{Time}_\textsc{Rttls}$ & $\mathrm{Time}_\textsc{Arttls}$
 & $\mathrm{Err}_\textsc{Rttls}$ & $\mathrm{Err}_\textsc{Arttls}$ \\ \hline
  Baart         & 4 & 8E-1   & 4  & 0.4970  & 0.0266   & 0.0567  & 2.22E-2  &4.24E-2\\
  Deriv2        & 7 & 2E-2   & 8  & 0.4062  & 0.0223   & 0.0335  & 4.39E-2  &1.45E-1\\
  Foxgood       & 3 & 5E-1   & 3  & 0.3086  & 0.0206   & 0.0302  & 3.48E-4  &2.14E-3\\
  Gravity       & 8 & 7E-1   & 8  & 0.3449  & 0.0246   & 0.0436  & 5.66E-3  &9.82E-3\\
  Heat          & 9 & 4E-1   & 8  & 0.3412  & 0.0275   & 0.0383  & 2.08E-1  &7.03E-2\\
  I\_Laplace    & 9 & 7E-1   & 10 & 0.6054  & 0.0408   & 0.1226  & 5.38E-3  &7.07E-2\\
  Phillips      & 7 & 4E-0   & 8  & 0.3218  & 0.0213   & 0.0367  & 1.44E-3  &7.80E-3\\
  Shaw          & 7 & 6E-1   & 7  & 0.6624  & 0.0330   & 0.0487  & 3.40E-3  &1.03E-2\\
           \hline
\end{tabular}
\caption{\label{tab:arttls} Algorithm \textsc{Arttls} on
ill-conditioned cases. The relative errors
$\mathrm{Err}_\textsc{Rttls} = || x_\textsc{rttls} - x_\textsc{ttls}
||_\infty / \|x_\textsc{ttls}\|_\infty$,
$\mathrm{Err}_\textsc{Arttls} = \| x_\textsc{ttls} -
x_\textsc{arttls} \|_\infty / \|x_\textsc{ttls}\|_\infty$. Both
\textsc{Rttls} and \textsc{Arttls} need less computational time than
\textsc{Ptls} based on Lanczos procedure.
}
\end{table}

\section{Conclusion}
In this paper, we derive a new perturbation bound for the total
least squares problem. This sharper and numerically computable
perturbation bound is well illustrated by the numerical examples.
 Also we show that three kinds of condition numbers  in \cite{Baboulin, Jia, Zhou}
obtained through different ways are mathematically equivalent.
%
We propose  randomized algorithms \textsc{Rtls}, \textsc{Rttls} and \textsc{Arttls}
for the numerical solutions of well-conditioned and ill-conditioned
total least squares problems, respectively. These randomized
algorithms can greatly reduce the computational time, and still give
solutions with good accuracy. The regularization parameter in
\textsc{Rttls} is estimated by the truncated parameter of the
\textsc{Tsvd} solution of $A x \approx b$ based on a fast randomized
\textsc{Svd} of $A$ \cite{Xiang}. Then a randomized \textsc{Svd} of
$[A,b]$ together with this truncation parameter yields a good
approximate \textsc{Ttls} solutions to the large-scale
ill-conditioned total least squares problems.
The detailed investigation on other regularization parameter
choices, and other techniques such as Tikhonov regularization, will
be our future research.

\section{Acknowledgments} The authors would like to thank Marc Baboulin, Ken Hayami, Zhongxiao Jia,
Lothar Reichel, and Jun Zou for their
comments and suggestions which led to improvements of our manuscript.

{\small
\bibliographystyle{siam}

}
\end{document}